\documentclass[twoside]{article}
\usepackage[usenames,dvipsnames]{color}
\usepackage{graphicx}
\usepackage{camnum}
\usepackage{harvard}

\usepackage{amsmath}
\usepackage{eufrak}
    \usepackage{ulem}
\newtheorem{remark}{Remark}

\usepackage{tikz,pgfplots}
        \pgfplotsset{compat = 1.3}
        \pgfplotsset{minor grid style={dotted}} \pgfplotsset{major grid
        style={dashed}}

        \pgfplotsset{every x tick label/.append style={font=\footnotesize,
        yshift=0.25ex}}

        \pgfplotsset{every y tick label/.append
        style={font=\footnotesize, xshift=0.25ex}}

        \definecolor{colorclassyorange}{rgb}{0.95000,0.32500,0.09800}
        \definecolor{colorchromeyellow}{rgb}{1.00000,0.6549,0}%
        \definecolor{colorpaleyellow}{rgb}{1.00000,0.8549,0.1}%

        \definecolor{colorclassyblue}{rgb}{0.00000,0.44706,0.74118}%
        \definecolor{colorpurple}{rgb}{0.49400,0.18400,0.55600}%
        \definecolor{colorfuschia}{rgb}{0.95039,0.0,0.95039}%
        \definecolor{colorlemongreen}{rgb}{0.6,0.8,0}%
        \definecolor{colorreal}{rgb}{0.92941,0.79412,0.12549}%
        \definecolor{colorimag}{rgb}{0.00000,0.49804,0.00000}%
        \definecolor{colorabs}{rgb}{1.00000,0.00000,0.00000}%

\newcommand{\ee}{{\mathrm e}}
\newcommand{\ii}{{\mathrm i}}

\newcommand{\K}{\mathcal{K}}

\newcommand{\FFF}{\mathcal{F}}
\newcommand{\AAA}{\mathcal{A}}
\newcommand{\DDD}{\mathcal{D}}

\newcommand{\dx}{\partial_x}

\def\O#1{{\cal O}\left(#1\right)}

\newcommand{\Int}[4]{\int_{#2}^{#3}\!#4\,\mathrm{d}#1}
\newcommand{\Ang}[2]{\left\langle #2 \right\rangle_{#1}}

\newcommand{\ang}[2]{\left\langle #1 \right\rangle_{#2}}

\newcommand{\schr}{Schr\"odinger }

\newcommand{\Hc}[1]{#1}

\newcommand{\dt}{\Hc{h}}

\newcommand{\mTh}[1]{\Hc{\Theta_{#1}(\dt)}}

\newcommand{\mom}[2]{\Hc{\mu_{#1,#2}(\dt)}}
\newcommand{\momO}{\Hc{\mu_{0,0}(\dt)}}
\newcommand{\Lf}[3]{\Hc{\Lambda\left[#1\right]_{#2,#3}\!(\dt)}}

\newcommand{\norm}[2]{\left\| #2 \right\|_{#1}}

\numberwithin{equation}{section}

\begin{document}
\title{Magnus--Lanczos methods with simplified commutators for the \schr equation with a time-dependent potential}

\author{Arieh Iserles,\footnote{Department of Applied Mathematics and Theoretical Physics,
 University of Cambridge, Wilberforce Rd, Cambridge CB3 0WA, UK.}\ \
  Karolina Kropielnicka\footnote{Institute of Mathematics, Polish Academy of Sciences,
  8 \'Sniadeckich Street, 00-656 Warsaw, Poland.}\ \  \&
  Pranav Singh\footnote{Mathematical Institute, Andrew Wiles Building, University of Oxford, Radcliffe Observatory Quarter, Woodstock Rd, Oxford OX2 6GG, UK and Trinity College, University of Oxford, Broad Street, Oxford OX1 3BH, UK.}}
\maketitle

\vspace*{6pt}

\vspace*{6pt} \noindent   \textbf{Keywords:} Schr\"odinger equation,
time-dependent potential, Magnus expansion, simplified commutators,
integral-preserving, Lanczos iterations, anti-commutators, Lie algebra,
oscillatory potentials, large time steps

\begin{abstract}
    The computation of the \schr equation featuring time-dependent
    potentials is of great importance in quantum control of atomic and
    molecular processes. These applications often involve highly
    oscillatory potentials and require inexpensive but accurate solutions
    over large spatio-temporal windows. In this work we develop
    Magnus expansions where commutators have been simplified. Consequently,
    the exponentiation of these Magnus expansions via Lanczos
    iterations is significantly cheaper than that for traditional Magnus
    expansions. At the same time, and unlike most competing methods, we
    simplify integrals instead of discretising them via quadrature at the
    outset -- this gives us the flexibility to handle a variety of
    potentials, being particularly effective in the case of highly
    oscillatory potentials, where this strategy allows us to consider
    significantly larger time steps.
%
%
%
\end{abstract}

\newpage

\setcounter{section}{0}

\section{Introduction}
We consider the linear, time-dependent \schr equation (TDSE) for a single particle moving in a time-varying electric field,
\begin{equation}\label{eq:schr}
\frac{\partial u(x,t)}{\partial t}=\ii \frac{\partial^2 u(x,t)}{\partial x^2} - \ii V(x,t) u(x,t),\quad x\in \BB{R},\ t \geq 0,
\end{equation}
where the complex-valued wave function $u=u(x,t)$ is given with an initial
condition $u(x,0) = u_0(x)$. Here $V(x,t)$ is a real-valued, time-dependent
electric field, and we are working in atomic units, where Planck's
constant is scaled to one ($\hbar =1$).

These equations are of great practical importance since they allow us to
study the behaviour of particles under the influence of changing electrical
field. As our ability to manipulate electric fields becomes more refined,
including the shaping of laser pulses, unprecedented quantum control of
atomic and molecular systems is becoming possible \cite{shapiro03qc}. Optimal
control of quantum systems is among the many challenges that require highly
accurate and computationally  inexpensive solutions of this equation, often
involving highly oscillatory potentials over large spatio-temporal windows

\subsection{Existing approaches}
\label{sec:existing} Time-dependent potentials significantly complicate
matters insofar as numerical solutions are concerned. Typically the solution
of \R{eq:schr} involves a truncation of the Magnus expansion, which is an
infinite series of nested integrals of nested commutators, as we will see in
this section.

Traditional methods for solving \R{eq:schr} usually commence with spatial discretisation,
\begin{equation}\label{eq:SD}
 \MM{u}'(t)=\ii (\K_2 -  \DDD_{V(\cdot,t)})\MM{u}(t),\qquad t\geq0,
\end{equation}
where the vector $\MM{u}(t)\in\BB{C}^M$ represents an approximation to the solution at time $t$, $\MM{u}(0)= \MM{u}_0$ is derived from the initial conditions, while $\K_2$ and $\DDD_{V(\cdot,t)}$ are $M \times M$ matrices which represent (discretisation of) second derivative and a multiplication by the interaction potential $V(\cdot,t)$, respectively.

{\bf Magnus expansions}. The system of ODEs (\ref{eq:SD}), which is of the
form
\begin{equation}\label{eq:Amag}
 \MM{u}'(t)= A(t) \MM{u}(t),\qquad t\geq0,
\end{equation}
with $A(t) = \ii (\K_2 -  \DDD_{V(\cdot,t)})$, can be solved via the Magnus
expansion \cite{magnus54},
\begin{equation}\label{eq:Amag2}
 \MM{u}(t)=\ee^{{\Theta}(t,s)}\MM{u}(s),
\end{equation}
where  ${\Theta}(t,s)$ is a time-dependent $M \times M$ skew-Hermitian matrix whose exponential evolves the solution from time $s$ to $t$. The Magnus expansion $\Theta(t,s)$ is obtained as an infinite series $\sum_{k=1}^\infty {\Theta}^{[k]}(t,s)$ with each ${\Theta}^{[k]}(t,s)$ composed of $k$ nested integrals and $k-1$ nested commutators (see expressions below).

In practice, we work with finite truncations of the Magnus series,
$$\Theta_m(t,s) = \sum_{k=1}^m {\Theta}^{[k]}(t,s),$$
and propagate the solution in suitably small time steps $\dt$,
$$\MM{u}^{n+1} = \ee^{\Theta_m(t_n+\dt,t_n)} \MM{u}^n,$$
in order to keep the truncation error low.
For the sake of simplicity, we analyse only  the first step,
\begin{equation}\label{eq:mag_disc}
\MM{u}^{1} = \ee^{\Theta_m(\dt)} \MM{u}^0,
\end{equation}
writing $\Theta_m(\dt)=\Theta_m(\dt,0)$, for short\footnote{The corresponding solution for any step $\ee^{\Theta_m(t_n+\dt,t_n)}  \MM{u}^n$ can be obtained by replacing $A(\zeta)$ by $A(t_n+\zeta)$ in the Magnus expansion.}. The first few terms of $\Theta_m(\dt)$ are
\begin{Eqnarray*}
  \Theta^{[1]}(\dt)=& &\int_0^\dt A(\xi_1)d\xi_1,\\
  \Theta^{[2]}(\dt)=& &-\frac{1}{2}\int_0^\dt \left[\int_0^{\xi_1} A(\xi_2)d\xi_2,A(\xi_1)\right]d\xi_1,\\
  \Theta^{[3]}(\dt)=& &\frac{1}{12}\int_0^\dt\left[\int_0^{\xi_1}A(\xi_2)d\xi_2,\left[\int_0^{\xi_1} A(\xi_2)d\xi_2,A(\xi_1)\right]\right]d\xi_1\\
 &&\mbox{} +\frac{1}{4}\int_0^\dt\left[\int_0^{\xi_1}\left[\int_0^{\xi_2}A(\xi_3)d\xi_3,A(\xi_2)\right]d\xi_2, A(\xi_1)\right]d\xi_1.\\
\end{Eqnarray*}

{\bf Exponential midpoint rule}. The simplest method in this family results
from letting $\Theta_1(h) = \Theta^{[1]}(\dt)$,
$$\MM{u}^{1} = \ee^{\Theta_1(h)} \MM{u}^0 = \exp\!\left(\int_0^\dt A(\xi) d\xi \right) \MM{u}^0 = \exp\!\left(\ii  \dt \K_2 - \ii  \int_0^\dt \DDD_{V(\xi)} d\xi \right) \MM{u}^0.$$
This method, called the exponential midpoint rule, is well known and has been
used for a long while \cite{TalEzer1992,lubich08fqc}. It is typical to
approximate $\int_0^\dt \DDD_{V(\xi)} d\xi = \DDD_{\int_0^\dt V(\xi) d\xi}$
by taking the value of $V$ at the middle of the integral, $\int_0^\dt V(\xi)
d\xi \approx \dt V(\dt/2)$, and concluding with an application of the Strang
splitting,
\begin{equation}
 \label{eq:expmidpoint}
 \MM{u}^1 = \exp\!\left(\Frac12 \ii \dt \K_2\right) \exp\!\left(- \ii \dt \DDD_{V(\dt/2)} \right) \exp\!\left(\Frac12 \ii \dt \K_2\right) \MM{u}^0.
\end{equation}

\begin{remark}\label{rmk:oddpowers}
The {\em power-truncated} Magnus expansions used here are time-symmetric and can be expanded solely in odd powers of $h$ \cite{iserles99ots,iserles00lgm}.
\end{remark}

\begin{remark}\label{rmk:quadrature}
We note that the first truncation, $\Theta_1$,
carries an error of $\O{\dt^3}$ \cite{iserles00lgm}. Ideally it should be
combined with an $\O{\dt^3}$ quadrature. However, due to Remark~\ref{rmk:oddpowers}, an $\O{\dt^{2n}}$ accuracy quadrature method will
have an accuracy of $\O{\dt^{2n+1}}$ in this context, resulting in the need for fewer quadrature points.
\end{remark}

As a consequence of Remark~\ref{rmk:quadrature} and the $\O{\dt^3}$ accuracy of the Strang splitting, \R{eq:expmidpoint} has a local error of $\O{\dt^3}$. %
%

{\bf Higher order truncations of the Magnus expansion}. Once higher order
accuracy is desired \cite{TalEzer1992,karlsson}, we need to consider higher
order truncations of the Magnus expansion such as
\[  \Theta_2(\dt) = \int_0^\dt A(\xi_1)d\xi_1 -\frac{1}{2}\int_0^\dt \left[\int_0^{\xi_1} A(\xi_2)d\xi_2,A(\xi_1)\right]d\xi_1 = \Theta(\dt) +\O{\dt^5}. \]
Higher order truncations necessarily involve nested integrals of nested commutators.
The nested integrals here need to be approximated using quadrature formul\ae{} of accuracy $\O{\dt^5}$. However, as mentioned before, it suffices to consider the Gauss--Legendre quadrature at only two nodes: $\tau_k = \Frac{\dt}{2}(1 \pm 1/\sqrt{3})$. This results in the method
\[ \MM{u}^{1} = \exp\!\left(\Frac{\dt}{2} (A(\tau_{-1})+A(\tau_{1})) + \Frac{\sqrt{3} \dt^2}{12} [A(\tau_{-1}), A(\tau_{1})] \right) \MM{u}^0. \]
 For the \schr equation \R{eq:schr}, this translates to

 \[ \MM{u}^{1} = \exp\!\left(\ii \dt \K_2 - \ii \dt \overline{V} + \dt^2 \Frac{\sqrt{3} \dt^2}{12} [\K_2, \widetilde{V}] \right) \MM{u}^0. \]
where $\overline{V} = [V(\tau_{-1})+V(\tau_{1})]/2$ and $\widetilde{V}= V(\tau_{-1})-V(\tau_{1})$.

{\bf Splitting the exponential of Magnus expansions}. The exponential of
$\Theta_2$ needs to be evaluated up to an accuracy of $\O{\dt^5}$. The
second-order Strang splitting,
\[ \ee^{\Theta_2(\dt)} = \ee^{\frac12 \ii \dt \K_2}\ \ee^{- \frac12 \ii \dt \DDD_{\overline{V}}}\ \ee^{\dt^2 \frac{\sqrt{3} \dt^2}{12} [\K_2,\DDD_{\widetilde{V}}]}\ \ee^{- \frac12 \ii \dt \DDD_{\overline{V}}}\ \ee^{\frac12 \ii \dt \K_2} + \O{\dt^3},\]
therefore, does not suffice.
Instead we require the fourth order Yoshida splitting, obtained by composing three order-two Strang splittings.

When the exponent to be split consists of two components, the number of
exponentials in an order-$2m$ Yoshida splitting grows as $2 \times 3^{m-1} +
1$. Here, we need to approximate the exponential of higher order Magnus
truncations, $\Theta_m$, which feature an increasingly larger number of
terms. Consequently the number of exponentials in the Yoshida splitting for
Magnus expansions grows even more rapidly.

Moreoever, we are left with the problem of evaluating the exponential of
commutators such as $[\K_2,\DDD_{\widetilde{V}}]$ which are expensive to
compute and do not posses a structure that allows for an easy exponentiation.
In higher-order methods such as $\Theta_3$ we start encountering commutators
in a nested form.

{\bf Magnus--Lanczos schemes.} An alternative approach for approximating the
exponential of the Magnus expansion is via Lanczos iterations
\cite{saad92kry}, leading to the popular Magnus--Lanczos schemes. This is,
arguably, a more flexible approach since we only require a method for
computing matrix--vector products of the form $\Theta_n v$ in each Lanczos
iteration.

Nevertheless the exponential growth resulting from the presence of nested
commutators is inevitable.
%
%
Moreover, the highly promising superlinear decay of error in the case of the
Lanczos method for approximating the matrix exponential is not seen until the
number of iterations is larger than the spectral radius of $\Theta_m(\dt)$
\cite{hochbruck97ksa}, which is very large unless the time step $\dt$ is
suitably small.

{\bf Commutator-free, integral-free quasi-Magnus methods}. To avoid the
exponential growth of cost due to presence of nested commutators, many
attempts have been made at deriving commutator-free schemes.
These usually proceed by replacing nested integrals in the
Magnus expansion by some quadratures or Taylor expansions of $V$ at the
outset, subsequently seeking a commutator-free exponential splitting that
adequately approximates the exponential of the discretised Magnus expansion.
Since the Magnus expansion does not appear explicitly in these schemes, they
are sometimes also referred to as {\em quasi-Magnus}.

For example in \cite{alvermann2011hocm}, instead of the exponential of Magnus
expansion, authors derive an alternative numerical propagator for \schr
equations, namely a product of exponentials of linear combinations of various
values of Hamiltonian operator (more precisely, values of Hamiltonian
operator are taken in Gauss-Legendre quadrature points).
\citeasnoun{blanes17quasimagnus}, on the other hand, investigate the
commutator-free expansion for differential equations of both parabolic and
hyperbolic equations, also providing stability and error analysis.

{\bf Integral-free Magnus--Zassenhaus splittings with simplified commutators}. In
\cite{BIKS} an integral-free numerical integrator with simplified commutators was proposed
for  \schr equation in the semiclassical regime. Once again, these proceed
via discretisation of the integrals in the Magnus expansion. However, unlike the commutator-free methods
where commutators are eliminated, these work by simplifying
the commutators in the Lie algebra of anti-commutators, subsequently
exploting the idea of the symmetric Zassenhaus asymptotic splittings
\cite{bader14eaf} for exponentiation. These have been shown to be highly
effective in the semiclassical regime.

{\bf Other notable approaches.} Expansion in Chebyshev polynomials is an
effective alternative to Lanczos iterations for approximating the exponential
of the Magnus expansion, particularly when large time steps are involved
\cite{TalEzer1992,NdongEzerKosloff10,TalEzer2012,Schaefer}. Other competing
approaches that forego Magnus expansions entirely include polynomial
approximations to the propagator based on Taylor expansions
\cite{LauvergnatTaylor07}, symplectic splitting methods
\cite{BlanesCasasMurua17}, Runge--Kutta methods \cite{TremblayCarrington04},
symplectic partitioned Runge--Kutta methods \cite{SanzSerna96} and the $(t,t')$
method \cite{Peskin}.

%
%

\subsection{Main contributions}
\label{sec:novelty}
In this work we present {\it Magnus--Lanczos methods with simplified commutators} that
\begin{enumerate}
\item
retain all the advantages of Magnus expansions and Lanczos methods,
\item
are free of nested commutators (and the associated growth in cost),
\item feature non-nested anti-commutators which preserve skew-Hermiticity
    of the expansion (thus unitarity of solution and stability of the
    discretised method),
\item preserve the integrals intact until the very last moment of the
    algorithm (this enables more flexibility, higher accuracy and often
    lower cost while dealing with numerical integration),
\item feature fewer nested integrals (due to identities \R{eq:II1} and
    \R{eq:II2}). Our order six methods, for instance, feature only
    twice-nested integrals instead of the quadruply-nested integrals that
    feature in a standard Magnus expansion -- not only is this beneficial
    for numerical quadrature, but it also makes analytic and asymptotic
    approximation easier.
\end{enumerate}

As discussed previously, standard Magnus expansions feature nested commutators. When we need to compute $\Theta_m v$ in each Lanczos iteration, these nested commutators result in the cost of  $\Theta_m v$
growing exponentially in $m$. The absence of nested commutators in our Magnus expansions, where such commutators have been simplified, results in the cost of $\Theta_m v$ growing linearly in $m$.

Moreover, we are able to do this while keeping integrals intact, resulting in
methods that are highly flexible -- not only it is possible to approximate
the integrals through any quadrature method, but we may also use exact
integrals for potentials whenever possible. This proves particularly
effective in the case of potentials with high temporal oscillations where we
no longer require a severe depression of time steps.



%

\subsection{Organisation of the paper}
Section 2 is devoted to the simplification of commutators in the Magnus expansion.
Magnus expansions for the \schr equation  evolve in the Lie algebra generated
by $\dx^2$ and $V(\cdot)$. However, as it will be pointed out in Subsection
2.1, a straightforward simplification of commutators of these operators using
the chain rule results in the loss of unitarity of the solution upon
discretisation. One of the novelties of our approach is working in the
algebra of anti-commutators, which leads to the preservation of
skew-Hermitian structure and stability of the scheme. The procedure for
deriving Magnus expansions with simplified commutators is presented in Subsections
2.2--2.6. In Subsections 2.7 and 2.8 we present concrete order four and order
six Magnus expansions, $\Theta_2$ and $\Theta_4$, respectively (methods
\R{eq:Omt3} and \R{eq:Th4}).

Section 3 deals with the implementation of our method. In
Subsection~\ref{sec:disc} we provide some details concerning spatial
discretisation. Subsection~\ref{sec:quadrature} deals with the evaluation of
derivatives and integrals of the potential appearing in the Magnus expansion.
While various alternatives are possible at this stage, a particular option --
namely, finite differences for derivatives and Gauss--Legendre quadrature
with three knots for integrals -- is outlined in greater detail (expressions
\ref{eq:FDu1}--\ref{eq:Lf}). In Subsection~\ref{sec:exponential} we discuss
the implementation of Lanczos iterations (achieved via \R{eq:Th4v}) for
numerical exponentiation of the Magnus expansion and the number of Fast
Fourier Transforms (FFTs) required per iteration.

Numerical examples are provided in Section 4, while in the last section we
briefly summarise our results.

\label{sec:layout}

\section{Magnus expansions with simplified commutators}
In contrast to the traditional approach of resorting to spatial
discretisation of \R{eq:schr}, which leads to the system of finite
dimensional ODEs $\R{eq:SD}$ followed by the Magnus expansion \R{eq:Amag2},
we begin straight away with a Magnus expansion of \R{eq:schr} while keeping
the underlying operators intact.

In numerically solving \R{eq:schr}, consistently with standard practice we
 impose periodic boundary conditions on a finite interval $I
\subseteq \BB{R}$. We further assume throughout that the interaction
potential $V(\cdot,t)$ and the wavefunction $u(\cdot,t)$ are sufficiently
smooth. For the purpose of this paper and for simplicity sake we assume that
they belong to $\CC{C}_p^\infty(I; \BB{R})$ and $\CC{C}_p^\infty(I; \BB{C})$,
respectively, the spaces of real valued and complex valued smooth periodic
functions over $I$, but our results extend in a straightforward manner to
functions of lower smoothness (the regularity constraints will depend on the
desired order of the method being derived).

Considering (\ref{eq:schr}) as an evolutionary PDE evolving in a Hilbert space, say $\mathcal{H} = \CC{L}^2(I;\BB{C})$, and suppressing the dependence on $x$,
\begin{equation}\label{eq:schrN}
\partial_t u(t) = \left(\ii \dx^2 - \ii  V(t)\right) u(t),\quad u(0)=u_0 \in \mathcal{H},
\end{equation}
is seen to be of the `ODE-like' form
\begin{equation}\label{eq:ODE}
\partial_t u(t) = A(t) u(t),\quad u(0)=u_0 \in \mathcal{H},
\end{equation}
with $A(t) = \ii \dx^2 - \ii  V(t)$. The operator $A(t)$ belongs to $\mathfrak{u}\!\left(\mathcal{H}\right)$, the Lie algebra of (infinite-dimensional) skew-Hermitian operators acting on the Hilbert space $\mathcal{H}$. Its flow is, therefore, unitary and resides in $\mathcal{U}\!\left(\mathcal{H}\right)$ -- the Lie group of unitary operators.

Unitary evolution of the wave function $u(t)$ under this flow is fundamental to quantum mechanics. Preservation of this property under discretisation is very important and we seek appropriate geometric numerical integrators to guarantee it. This comes about naturally once we work in the correct Lie-algebraic framework. As we note later, unitarity also guarantees stability of a consistent numerical scheme.


For a general equation of the form (\ref{eq:ODE}) where $A(t)$ resides in a Lie algebra $\mathfrak{g}$, the solution for the flow can be formally written in the form of a Magnus expansion,
\begin{equation}\label{eq:magnus_exp_t}
u(\dt) = \ee^{\Theta(\dt)}u(0),
\end{equation}
which differs from \R{eq:Amag2} in the sense that the Magnus expansion $\Theta$ is in general an infinite-dimensional and unbounded operator, not a matrix.

\begin{remark}
Convergence of the Magnus expansion, in the finite dimensional case is only
guaranteed for sufficiently small time steps \cite{moan08magnus}. In
principle, this becomes problematic when we consider Magnus expansions of
undiscretised and unbounded operators. An extension to bounded operators on an
infinite-dimensional Hilbert space was done by \cite{casasMagnus}.

Rigorous analysis in the context of the \schr equation, which features an unbounded operator,
 has
been carried out by \cite{Hochbruck02onmagnus} who show that when Magnus
expansions of unbounded operators are considered in a formal sense, the
concrete methods based on these approaches do demonstrate the expected order
of convergence.

Pursuing a similar strategy, it is possible to derive rigorous error bounds
for the Magnus expansion based methods presented in this manuscript. However,
since this analysis involves development of additional theory that could
obscure the presentation of the proposed methods, it will be beyond the scope
of our investigations. Here we refer the curious reader to Chapter 9. of
\cite{Singh17}. \end{remark}

\subsection{The algebra of anti-commutators}
The vector field in the \schr equation \R{eq:schrN} is a linear combination of the action of two operators, $\dx^2$ and the multiplication by the interaction potential $V(t)$, for any $t \geq 0$. Since the Magnus expansion requires nested commutation, the focus of our interest is the Lie algebra generated by $\dx^2$ and $V(\cdot)$,
\begin{displaymath}
  \GG{F}=\CC{LA} \{\dx^2, V(\cdot) \},
\end{displaymath}
i.e.\ the  linear-space closure of all nested commutators of $\dx^2$ and $V(\cdot)$.

{\bf Simplifying commutators}. To simplify commutators, we follow the
approach of \cite{bader14eaf} and study their action on functions. For
example, using the chain rule we find,
\begin{displaymath}
  [\dx^2,V]u=\dx^2(Vu) - V(\dx^2u)=(\dx^2V)u+2(\dx V)\dx u,
\end{displaymath}
which implies that $[\dx^2,V]=(\dx^2V)+2(\dx V)\dx$. Similarly, we conclude that
\begin{align*}
  [\dx^2,[\dx^2,V]]&=(\dx^4V)+4(\dx^3V)\dx+4(\dx^2V)\dx^2,\\
  [V,[\dx^2,V]] &= -2(\dx V)^2.
\end{align*}
Note that we have ignored here the dependence on the time variable since the derivatives are only in the spatial variable.

{\bf Loss of skew-Hermiticity}. Simplifying commutators in this way, we can,
in principle, get rid of all nested commutators occurring in the (truncated)
Magnus expansion of the undiscretised operators. It is only after this stage
that we would resort to spatial discretisation. Proceeding in this way,
however, we break an important structural property -- upon discretisation,
such a Magnus expansion is no longer skew-Hermitian and thus
its exponential is no longer unitary.

To illustrate the loss of unitarity, let us consider the two differential operators: $[\dx^2,V]$ and $(\dx^2V)+2(\dx V)\dx$.  Spatial discretisation transforms these two analytically identical operators to $[\K_2, \DDD_V ]$ and $ \DDD_{\dx^2 V} + 2 \DDD_{\dx V} \K_1$, respectively.
Assuming that $\K_n$ is skew-symmetric for odd $n$ and symmetric for even $n$ (i.e. the skew-symmetry of $\dx$ is preserved under discretisation) and $\DDD_V$ is symmetric (note that $V$ is real-valued and $\DDD_V$ represents multiplication by $V$), the commutator $[\K_2,\DDD_V ]$ is skew-symmetric. However, the second expression, subject to discretisation, is no longer skew-symmetric.
A similar problem is encountered in the discretisation of $[\dx^2,[\dx^2,V]]$ following the simplification of the commutator.

The loss of skew-symmetry in the simplification of $[\dx^2,V]$ (and
skew-Hermiticity in general) is a cause for concern on two accounts: firstly,
the exponential of a Magnus expansion which features terms like $[\dx^2,V]$
($\Theta_2$, for instance) is no longer unitary, which is highly undesirable
insofar as the physics is concerned; secondly, since $ \DDD_{\dx^2 V} + 2
\DDD_{\dx V} \K_1 $ has (large) real eigenvalues, its exponential blows up,
which is highly undesirable from the numerical point of view. This blowup can
be extreme even in the simplest of cases (see Figure~\ref{fig:stability}).

\begin{figure}[tbh]
\begin{center}
%
%
\begin{tikzpicture}

\begin{axis}[%
width=2in,
height=1.5in,
scale only axis,
scale only axis,
xmode=log,
xmin=1e-06,
xmax=1,
xminorticks=true,
ymode=log,
ymin=1,
ymax=10000,
xlabel={$t$},
ylabel={$\norm{2}{\exp(t \AAA)}$},
yticklabel style = {font=\small,xshift=0.5ex},
xticklabel style = {font=\small,yshift=0.5ex},
yminorticks=true,
legend style={legend cell align=left,align=right,draw=black,font=\footnotesize},
legend pos=north west
]
\addplot [color=black,line width=1.0pt,solid,mark=o,mark size=1,mark options={solid},mark repeat=5]
  table[row sep=crcr]{%
9.5367431640625e-07	1.00000000000001\\
1.34869915234861e-06	1.00000000000001\\
1.9073486328125e-06	1\\
2.69739830469722e-06	1.00000000000001\\
3.814697265625e-06	1.00000000000001\\
5.39479660939444e-06	1.00000000000001\\
7.62939453125e-06	1\\
1.07895932187889e-05	1.00000000000001\\
1.52587890625e-05	1\\
2.15791864375777e-05	1.00000000000001\\
3.0517578125e-05	1\\
4.31583728751555e-05	1\\
6.103515625e-05	1\\
8.6316745750311e-05	1\\
0.0001220703125	1\\
0.000172633491500622	1\\
0.000244140625	1\\
0.000345266983001244	1\\
0.00048828125	1\\
0.000690533966002488	1.00000000000002\\
0.0009765625	1.00000000000003\\
0.00138106793200498	1.00000000000005\\
0.001953125	1.00000000000006\\
0.00276213586400995	1.00000000000009\\
0.00390625	1.00000000000012\\
0.0055242717280199	1.00000000000016\\
0.0078125	1.00000000000021\\
0.0110485434560398	1.00000000000029\\
0.015625	1.00000000000039\\
0.0220970869120796	1.00000000000053\\
0.03125	1.00000000000072\\
0.0441941738241592	1.00000000000099\\
0.0625	1.00000000000137\\
0.0883883476483184	1.00000000000188\\
0.125	1.00000000000261\\
0.176776695296637	1.00000000000359\\
0.25	1.00000000000495\\
0.353553390593274	1.00000000000678\\
0.5	1.00000000000932\\
0.707106781186547	1.00000000001279\\
1	1.00000000001752\\
};
\addlegendentry{$[\K_2,\DDD_V]$};

\addplot [color=colorclassyorange,solid,line width=1.0pt]
  table[row sep=crcr]{%
9.5367431640625e-07	1.00121992249111\\
1.34869915234861e-06	1.00172566592511\\
1.9073486328125e-06	1.00244132983919\\
2.69739830469722e-06	1.00345430027724\\
3.814697265625e-06	1.00488859287324\\
5.39479660939444e-06	1.00692045651732\\
7.62939453125e-06	1.00980086787215\\
1.07895932187889e-05	1.01388819176601\\
1.52587890625e-05	1.01969604639192\\
2.15791864375777e-05	1.02796428703858\\
3.0517578125e-05	1.03976579098561\\
4.31583728751555e-05	1.05666969267405\\
6.103515625e-05	1.08099481921115\\
8.6316745750311e-05	1.11620721498467\\
0.0001220703125	1.16754034057858\\
0.000172633491500622	1.24292044016087\\
0.000244140625	1.35415005237855\\
0.000345266983001244	1.51766514366917\\
0.00048828125	1.7521752104763\\
0.000690533966002488	2.06675110386046\\
0.0009765625	2.43458899874847\\
0.00138106793200498	2.78241525316906\\
0.001953125	3.08365353457636\\
0.00276213586400995	3.43826997091704\\
0.00390625	3.78559892649994\\
0.0055242717280199	4.16231860197933\\
0.0078125	4.56241984040098\\
0.0110485434560398	4.97326111978984\\
0.015625	5.43530658068583\\
0.0220970869120796	5.91111879597626\\
0.03125	6.42554211383186\\
0.0441941738241592	6.99992819856018\\
0.0625	7.57654210796593\\
0.0883883476483184	8.11679944211834\\
0.125	8.89221766834314\\
0.176776695296637	9.84873026675045\\
0.25	11.0883587453842\\
0.353553390593274	13.7386094877787\\
0.5	36.0000561856821\\
0.707106781186547	154.801099455226\\
1	1230.11721079546\\
};
\addlegendentry{$\DDD_{\dx^2 V}+ 2 \DDD_{\dx V} \K_1$};
%

\end{axis}
\end{tikzpicture}%
\end{center}
    \caption{The two equivalent forms $[\dx^2,V]$ and  $(\dx^2V) + 2(\dx V)\dx$,
    lead to two different discretisations, $[\K_2,\DDD_V]$ and $ \DDD_{\dx^2 V} + 2 \DDD_{\dx V} \K_1 $.
    While the exponential of the former is unitary, the exponential of the latter blows up.}
    \label{fig:stability}
\end{figure}

\begin{remark}
Note that $[\dx^2,V]$ and the simplified form $(\dx^2V)+2(\dx V)\dx$ are both
 skew-Hermitian operators. However, the straightforward discretisation (i.e.\
where every instance of $V^{(k)}$ is replaced by $\DDD_{V^{(k)}}$ and $\dx^k$ by $\K_k$)
of the simplified form to $\DDD_{\dx^2 V} + 2 \DDD_{\dx V} \K_1$ is where
skew-Hermiticity is lost. This is not entirely surprising since a
discretisation scheme, in general, might respect the skew-Hermiticity of only
a subset of skew-Hermitian operators.
\end{remark}


{\bf The Lie algebra of anti-commutators}. In the methods presented here we
 circumvent the problem of skew-Hermitian
discretisation by working with differential operators of the form
\begin{equation}\label{eq:angbra}
  \ang{f}{k} := \Frac12 \left( f \circ \dx^k + \dx^k \circ f \right), \quad  k \geq 0,\ f \in \CC{C}_p^{\infty}(I;\BB{R}),
\end{equation}
which are inherently symmetrised. These are anti-commutators of $f$ and
$\dx^k$. The action of this differential operator on $u$, for example, is
$$ \ang{f}{k} u = \Frac12 \left( f \dx^k u + \dx^k (f u) \right) $$
and the discretisation of this operator is
\begin{equation}\label{eq:angbradisc}
 \ang{f}{k} \leadsto \Frac12 \left( \DDD_f \K_k + \K_k \DDD_f \right).
\end{equation}

\begin{remark}
It is a very simple process to verify that $\ang{f}{k}$ is a skew-Hermitian
operator for odd $k$ and Hermitian operator for even $k$. This feature is
maintained under straightforward discretisation. This is the reason why the
choice of algebra of anti-commutators $\ang{\cdot}{k}$ seems to be optimal
for our purposes.
\end{remark}

The commutators of these operators can be solved using the following rules,
\begin{Eqnarray}
\label{eq:id}\left[ \Ang{1}{f}, \Ang{0}{g} \right] & =& \Ang{0}{f (\dx g)},\\
\nonumber\left[ \Ang{1}{f}, \Ang{1}{g} \right] & =& \Ang{1}{f (\dx g)-(\dx f) g},\\
\nonumber\left[ \Ang{2}{f}, \Ang{0}{g} \right] & =& 2 \Ang{1}{f (\dx g)},\\
\nonumber\left[ \Ang{2}{f}, \Ang{1}{g} \right] & =& \Ang{2}{2 f (\dx g) -
(\dx f) g} - \Frac12 \Ang{0}{2(\dx f) (\dx^2 g) + f (\dx^3
g)},\\
\nonumber\left[ \Ang{2}{f}, \Ang{2}{g} \right] & =& 2 \Ang{3}{f (\dx g)
- (\dx f) g} + \Ang{1}{2(\dx^2 f) (\dx g) - 2(\dx f) (\dx^2 g) + (\dx^3 f) g
-
f (\dx^3 g)},\\
\nonumber\left[ \Ang{3}{f}, \Ang{0}{g} \right] & =& 3 \Ang{2}{f (\dx g)} - \Frac12 \Ang{0}{3(\dx f) (\dx^2 g) + f (\dx^3 g)},\\
\nonumber\left[ \Ang{4}{f}, \Ang{0}{g} \right] & =& 4 \Ang{3}{f (\dx g)} -2
\Ang{1}{3(\dx f) (\dx^2 g) + f (\dx^3 g)}.
\end{Eqnarray}

\begin{remark} Note that each commutator is simplified to a linear combination of $\ang{\cdot}{k}$s where the indices $k$ are either
even or odd, but the two are never mixed. Upon discretisation these
individual terms will result in either Hermitian matrices or skew-Hermitian
matrices, but never a mixture of the two. Moreover, any term can be discarded
(say, due to a small size) without disturbing the symmetry.
\end{remark}

There is rich algebraic theory underlying these operators (including a
general formula for  \R{eq:id}) which  feature in a separate publication
\cite{psalgebra}.  In principle, however, the above rules can be verified by
application of the chain rule.

\begin{remark}
Note that, by definition \R{eq:angbra},
\begin{enumerate}
    \item these brackets are linear, so that $\Ang{2}{2 f (\dx g) - (\dx f) g} = 2 \Ang{2}{f (\dx g)} - \Ang{2}{(\dx f) g}$,
    \item $\ang{f}{0} = f$; and
    \item $\ang{1}{2} = \dx^2$.
\end{enumerate}
\end{remark}

With this new notation in place and using \R{eq:id}, we can now simplify commutators to anti-commutators,
    \begin{Eqnarray}
        \nonumber[\ii\dx^2, \ii V ] & =& -[ \ang{1}{2}, \ang{V}{0} ] = - 2 \ang{\dx V}{1},\\
        \nonumber[\ii V,[\ii\dx^2,\ii V]] & = & -\ii [ \ang{V}{0},[ \ang{1}{2},\ang{V}{0}]] = 2\ii \ang{(\dx V)^2}{0},\\
        \nonumber[\ii \dx^2,[\ii \dx^2,\ii V]] & = & -\ii [\ang{1}{2},[\ang{1}{2}, \ang{V}{0} ]] = \ii \ang{\dx^4 V}{0} - 4 \ii \ang{\dx^2 V}{2}.
    \end{Eqnarray}
Straightforward discretisations of these operators preserve the symmetries that are crucial for preserving unitarity.

\subsection{The expansion for the \schr equation}
Using the approach introduced in the previous section, commutators in the Magnus expansion, $\Theta(\dt) = \sum_{k=1}^\infty \Theta^{[k]}(\dt)$, can be expanded in terms of anti-commutators: $\Theta(h)=\sum_{k=0}^\infty\ii^{k+1}\ang{f_k}{k}$. Since $A(t) = \ii \dx^2 - \ii  V(t) = \ii \ang{1}{2} - \ii  \ang{V(t)}{0}$,
the first term of the Magnus expansion is the integral $\Theta^{[1]}(h) =  \Int{\zeta}{0}{\dt}{A(\zeta)}$,
\begin{equation}
\label{eq:MEP1} \Theta^{[1]}(h) =  \ii  \dt \ang{1}{2} - \ii  \Int{\zeta}{0}{\dt}{\ang{ V(\zeta) }{0}} =  \ii  \dt \ang{1}{2} - \ii  \ang{\Int{\zeta}{0}{\dt}{ V(\zeta) }}{0}.
\end{equation}
Note that the integrals here are in time, while differential operators are in space. Along with linearity of the brackets and integrals, this observation allows us to exchange brackets and integrals.

{\bf The first non-trivial term}. $\Theta^{[2]}(h)$, is simplified as
\begin{Eqnarray}
\nonumber
 \Theta^{[2]}(h) &=& -\Frac{1}{2} \Int{\zeta}{0}{\dt}{\Int{\xi}{0}{\zeta}{\left[A(\xi), A(\zeta)  \right]}}\\
\nonumber &=&-\Frac{1}{2}\Int{\zeta}{0}{\dt}{\Int{\xi}{0}{\zeta}{\left[\ii \ang{1}{2} - \ii  \ang{V\left(\xi\right)}{0},\ \ii \ang{1}{2} - \ii  \ang{V\left(\zeta\right)}{0}  \right]}}\\
\nonumber &=&-\Frac{1}{2}\Int{\zeta}{0}{\dt}{\Int{\xi}{0}{\zeta}{\left[\ang{1}{2}, \ang{V\left(\zeta\right)}{0}\right] + \left[\ang{V\left(\xi\right)}{0}, \ang{1}{2} \right]}}\\
\nonumber &=&-\left(  \Int{\zeta}{0}{\dt}{\zeta \ang{\dx V(\zeta)}{1} } -  \Int{\zeta}{0}{\dt}{ \Int{\xi}{0}{\zeta}{ \ang{\dx V(\xi)}{1} } }  \right)\\
\label{eq:MEP2} &=&-\ang{\Int{\zeta}{0}{\dt}{\zeta (\dx V(\zeta)) } -  \Int{\zeta}{0}{\dt}{ \Int{\xi}{0}{\zeta}{(\dx V(\xi)) } }  }{1}.
\end{Eqnarray}

{\bf Higher order terms}. Similarly, using \R{eq:id}, we can simplify higher
nested commutators in the Magnus expansion. For instance,
\[ \Theta^{[3,1]}(h) = \frac{1}{12}\int_0^\dt\left[\int_0^{\xi_1}A(\xi_2)d\xi_2,\left[\int_0^{\xi_1} A(\xi_2)d\xi_2,A(\xi_1)\right]\right]d\xi_1, \]
which occurs as a part of $\Theta^{[3]}(h)$, is simplified to
\begin{Eqnarray}
\Theta^{[3,1]}(h) & = &\ \ \Frac{1}{3}\ii  \ang{\Int{\zeta}{0}{\dt}{ \zeta^{2} (\dx^{2}V(\zeta)) } - \Int{\zeta}{0}{\dt}{ \zeta \Int{\xi}{0}{\zeta}{ (\dx^{2}V(\xi)) } }}{2}\nonumber \\
\nonumber &&\mbox{}+\Frac{1}{6} \ii  \ang{ \Int{\zeta}{0}{\dt}{ \zeta (\dx V(\zeta)) \Int{\xi}{0}{\zeta}{ (\dx V(\xi)) } } - \Int{\zeta}{0}{\dt}{ \left( \Int{\xi}{0}{\zeta}{ (\dx V(\xi)) }\right)^2 }}{0}\\
\label{eq:MEP3}&&\mbox{} - \Frac{1}{12} \ii \ang{ \Int{\zeta}{0}{\dt}{ \zeta^{2} (\dx^{4}V(\zeta)) } -  \Int{\zeta}{0}{\dt}{ \zeta \Int{\xi}{0}{\zeta}{ (\dx^{4}V(\xi)) } }}{0}.
\end{Eqnarray}

\subsection{Simplification of integrals}

After simplifying terms in the Magnus expansion we arrive at expressions such as (\ref{eq:MEP2}) and (\ref{eq:MEP3}), where each integral is of the form
\[ I_{ \mathcal{S},  {\tiny \MM{f}} }(\dt) = \int_{\mathcal{S}}  f(\MM{\xi})\,\mathrm{d}\MM{\xi}, \]
where $\MM{f}(\MM{\xi})=\prod_{j=1}^s f_j(\xi_j)$ for some function $f_j$, and $\mathcal{S}$ is an $s$-dimensional polytope of the special form,
\begin{equation}
\label{eq:polytope} \mathcal{S} = \{\MM{\xi} \in \BB{R}^s\ :\ \xi_1 \in [0,\dt],\quad \xi_l \in [0, \xi_{m_l}],\quad l=2,3,\ldots,s \},
\end{equation}
where $m_l \in \{1,2,\ldots, l-1\}$, $l=2,3, \ldots, s$. For details about the types of polytopes of integration appearing in the Magnus expansion see \cite{iserles00lgm}.


The special form of these polytopes and the integrands obtained after expanding the commutators, allows us to simplify the terms of the Magnus expansion further. Integration by parts leads us to the following identities:
\begin{Eqnarray}
\label{eq:II1} &&\int_0^hf_1(\xi_1)\left(\int_0^{\xi_1}f_2(\xi_2)d\xi_2\right)d\xi_1=\int_0^hf_2(\xi_1)\left(\int_{\xi_1}^hf_1(\xi_2)d\xi_2\right)d\xi_1,\\
\label{eq:II2}  &&\int_0^hf_1(\xi_1)d\xi_1\left(\int_0^{\xi_1}f_2(\xi_2)d\xi_2\right)\left(\int_0^{\xi_1}f_3(\xi_3)d\xi_3\right)d\xi_1\\
\nonumber&=&\int_0^h\left(\int_{\xi_3}^hf_1(\xi_1)d\xi_1\right)\left(f_2(\xi_3)\int_0^{\xi_3}f_3(\xi_2)d\xi_2+f_3(\xi_3)\int_0^{\xi_3}f_2(\xi_2)d\xi_2\right)d\xi_3.
\end{Eqnarray}
In our simplifications, (\ref{eq:MEP2}) and (\ref{eq:MEP3}), we have already encountered integrals over a triangle such as $\Int{\zeta}{0}{\dt}{ \Int{\xi}{0}{\zeta}{ (\dx V(\xi)) } }$ and $\Int{\zeta}{0}{\dt}{ \zeta \Int{\xi}{0}{\zeta}{ (\dx^{2}V(\xi)) } }$. We can reduce these to integrations over a line by applying the first identity with $f_1(\xi_1)=1$, $f_2(\xi_2)=\dx V(\xi_2)$ and $f_1(\xi_1)=\xi_1$, $f_2(\xi_2)=\dx V(\xi_2)$, respectively. Integration over the pyramid in $\Int{\zeta}{0}{\dt}{ \left( \Int{\xi}{0}{\zeta}{ (\dx V(\xi)) }\right)^2 }$ is similarly reduced using the second identity with $f_1(\xi_1)=1$, $f_2(\xi_2)=\dx V(\xi_2)$, $f_3(\xi_3)=\dx V(\xi_3)$.

\begin{remark}
The use of identities \R{eq:II1} and \R{eq:II2} is what allows us to reduce the complexity of the integrals in our Magnus expansions. In particular, our order-six method in Subsection 2.8 features integrals over a triangle instead of integrals over four-dimensional polytopes that are typical in the usual order-six Magnus expansions.
\end{remark}

\begin{remark}
Although it might be possible to develop general formalism for extending these observations to higher dimensional polytopes appearing in the Magnus expansion, the two identities presented here suffice for all results presented in our work. Deducing similarly useful identities for reduction of nested integrals in any specific high dimensional polytope should also be possible and would be a very helpful result.
\end{remark}

\subsection{A proposed Magnus expansion}
After simplification of commutators and applications of the integration identities (\ref{eq:II1}) and (\ref{eq:II2}), the Magnus expansion $\Theta_3$, for instance, reduces to the sum of the following terms,
\begin{Eqnarray}
\label{eq:tree1} \Theta^{[1]}(h) & = & \ii  \dt \dx^2 - \ii  \Int{\zeta}{0}{\dt}{ V(\zeta) },\\
\label{eq:tree2} \Theta^{[2]}(h) & = & -2 \Ang{1}{ \Int{\zeta}{0}{\dt}{
\left(\zeta - \Frac{\dt}{2}\right) (\dx V(\zeta)) } },
\end{Eqnarray}

\begin{Eqnarray}
                 \Theta^{[3,1]}(h) & = & -\Frac{1}{6}\ii  \Int{\zeta}{0}{\dt}{ \Int{\xi}{0}{\zeta}{ \left( 2 \dt-3 \zeta \right)   (\dx V(\zeta)) (\dx V(\xi))} }  \nonumber\\
\label{eq:tree3} && \mbox{} - \Frac{1}{6}\ii \Ang{2}{ \Int{\zeta}{0}{\dt}{
\left( \dt^{2}-3 \zeta^{2} \right) (\dx^{2}V(\zeta)) } },\\
                 \Theta^{[3,2]}(h) & = & \Frac{1}{2}\ii   \Int{\zeta}{0}{\dt}{ \Int{\xi}{0}{\zeta}{ \left( \zeta-2 \xi \right)   (\dx V(\zeta)) (\dx V(\xi))} } \nonumber\\
\label{eq:tree4} && \mbox{}  +\Frac{1}{2}\ii \Ang{2}{ \Int{\zeta}{0}{\dt}{ \left( \dt^{2}-4 \dt \zeta+3 \zeta^{2} \right) (\dx^{2}V(\zeta)) } },
\end{Eqnarray}
where $\Theta^{[3,2]}(h)$ refers to the second part of $\Theta^{[3]}$. Here and in the sequel we prefer to express $\Ang{0}{f}$ as $f$ and $\Ang{2}{1}$ as $\partial_x^2$ to avoid an excessively pedantic and longwinded notation.

\begin{remark}
The term $\Theta^{[2]}(h) = -2 \Ang{1}{\Int{\zeta}{0}{\dt}{ \left(\zeta-\Frac{\dt}{2}\right) (\dx V(\zeta)) }} $
might seem to be $\O{\dt^2}$ at first sight.
A closer look at the special form of the integrand, however, shows that the term is, in fact, $\O{\dt^3}$. To observe this, consider $V(\zeta)$ expanded about $0$, so that $V(\zeta) = V(0)+ \sum_{k=1}^\infty \zeta^k  V^{(k)}(0)/{k!} $. Note that the $\dt^2$ term $ \Int{\zeta}{0}{\dt}{ \left(\zeta-\Frac{\dt}{2}\right) (\dx V(0))}$ vanishes. This is consistent with Remark~\ref{rmk:oddpowers}. Similar care has to be exercised throughout the simplifications while analysing size. We refer the reader to \cite{iserles00lgm} for a more general analysis of such gains of powers of $\dt$, which occurs in specific terms of the Magnus expansion due to their structure.
\end{remark}
%

%
%

\subsection{Simplifying notation}
The algebraic workings become increasingly convoluted once we start dealing with larger nested commutators and integrals. Here it becomes helpful to introduce a notation for the integrals on the line,
\begin{equation}
\label{eq:mom}
\mom{j}{k} = \Int{\zeta}{0}{\dt}{\tilde{B}_{j}^k(\dt,\zeta) V(\zeta) },
\end{equation}
and integrals over a triangle,
\begin{equation}
\label{eq:Lf1}
\Lf{f}{a}{b} = \Int{\zeta}{0}{\dt}{\Int{\xi}{0}{\zeta}{f(\dt,\zeta,\xi) \left[ \dx^a V(\zeta)\right] \left[\dx^b V(\xi) \right] }},
\end{equation}
where $\tilde{B}$ is a rescaling of Bernoulli polynomials \cite{abram64spfunc,lehmer88ber},
\begin{equation*}
\label{eq:Bernoulli}
\tilde{B}_j(\dt,\zeta) = \dt^j B_j\left(\zeta/\dt \right).
\end{equation*}

\subsection{Fourth order Magnus expansions}
With this new notation in place,
the Magnus expansions $\Theta_2(h)$ and $\Theta_3(h)$ can be presented more concisely,
\begin{Eqnarray}
\label{eq:Omt3}
\Theta_2(h) & =&  \overbrace{\ii \dt  \dx^2 - \ii    \momO}^{\O{\dt}} - \overbrace{2\Ang{1}{\dx \mom{1}{1}}}^{\O{\dt^3}},\\
\label{eq:Omt4}
\Theta_3(h)  &= &\Theta_2(h) + \overbrace{\ii   \Lf{\psi}{1}{1} +2 \ii \Ang{2}{\dx^{2} \mom{2}{1}}}^{\O{\dt^4}},
\end{Eqnarray}
where
\[ \psi(\dt,\zeta,\xi) = \zeta - \xi - \Frac{\dt}{3}. \]
Due to Remark~\ref{rmk:oddpowers}, the error in these Magnus expansions can be expanded solely in odd powers of $h$ since they are based on the power-truncated Magnus expansions of \cite{iserles00lgm}. Consequently, both of these expansions are fourth-order.

\begin{remark}
Since the $j${\em th} rescaled Bernoulli polynomial scales as $\O{\dt^j}$, we expect $\mom{j}{k} = \O{\dt^{jk+1}}$ . Since the integral of the Bernoulli polynomials vanishes,
\begin{equation}
\label{eq:IntB}
 \Int{\zeta}{0}{\dt}{B_j(\dt,\zeta) }= 0,
\end{equation}
however, the term $\mom{j}{1}$ gains an extra power of $\dt$ and is $\O{\dt^{j+2}}$ .

In general, for a polynomial $p_n(\dt,\zeta,\xi)$ featuring {\em only} degree-$n$ terms in $\dt, \zeta$ and $\xi$, the linear (integral) functional (\ref{eq:Lf1}) is $\O{\dt^{n+2}}$. However, the integral of $\psi$ over the triangle vanishes,
\begin{equation}
\label{eq:IntPhi}
 \Int{\zeta}{0}{\dt}{\Int{\xi}{0}{\zeta}{\psi(\dt,\zeta,\xi)}}=0,
\end{equation}
lending an extra power of $\dt$ to terms featuring $\Lf{\psi}{a}{b}$.
\end{remark}

\subsection{Sixth order Magnus expansion}
Arbitrarily high order Magnus expansions with simplified commutators can be derived by following the procedure described in the preceding sections. The order six expansion, $\Theta_4(h)$, for instance, is \begin{Eqnarray}\label{Theta_4}
\nonumber \Theta_4(h) & = &  \overbrace{\ii \dt \dx^2 - \ii  \momO}^{\O{\dt}} -  \overbrace{2 \Ang{1}{\dx \mom{1}{1}}}^{\O{\dt^3}} + \overbrace{\ii  \Lf{\psi}{1}{1} + 2\ii \Ang{2}{\dx^2 \mom{2}{1}}}^{\O{\dt^4}} \\
\nonumber && +\overbrace{\Frac{1}{6} \Ang{1}{\Lf{\varphi_1}{1}{2} + \Lf{\varphi_2}{2}{1} }}^{\O{\dt^4}} + \overbrace{\Frac{1}{6} \Ang{1}{\Lf{\phi_1}{1}{2} + \Lf{\phi_2}{2}{1} }}^{\O{\dt^5}}\\
\label{eq:Omt5}&&+\overbrace{\Frac{4}{3}  \Ang{3}{\dx^3 \mom{3}{1}}}^{\O{\dt^5}} +\overbrace{\Frac{1}{4} \ii \dx^4 \mom{2}{1}}^{\O{\dt^4}} = \mTh{} + \O{\dt^7},
\end{Eqnarray}
where
\begin{Eqnarray}
\varphi_1(\dt,\zeta,\xi)& = & \dt^2 - 4 \dt \xi + 2 \zeta \xi, \\
\nonumber \varphi_2(\dt,\zeta,\xi)& = & (\dt-2\zeta)^2 - 2 \zeta \xi,\\
\nonumber \phi_1(\dt,\zeta,\xi) & = & \dt^2 - 6 \dt \zeta + 6 \dt \xi + 6 \zeta \xi + 3 \zeta^2 - 12 \xi^2,\\
\nonumber \phi_2(\dt,\zeta,\xi) & = & \dt^2 - 6 \dt \zeta + 6 \dt \xi - 6 \zeta \xi + 5 \zeta^2.
\end{Eqnarray}
\begin{remark}
Integrals of $\phi_j$ vanish over the triangle,
\begin{equation}
\label{eq:IntTheta}
 \Int{\zeta}{0}{\dt}{\Int{\xi}{0}{\zeta}{\phi_j(\dt,\zeta,\xi)}}=0,\quad j= 1,2,
\end{equation}
lending an extra power of $\dt$ to the functionals wherever $\phi_j$s appear. No similar observation about $\varphi_j$s can be made and they have been kept separate in (\ref{Theta_4}) from the $\Lambda$ terms featuring $\phi_j$
simply since they lead to terms of different orders. For purposes of
computing, however, they can be combined.
\end{remark}
The linearity of the brackets means that $\Theta_4(h)$ can be rewritten in
the form,
\begin{equation}
\label{eq:Th4}
\Theta_4(h) = \theta_0 + \ang{\theta_1}{1} + \ang{\theta_2}{2} + \ang{\theta_3}{3}
\end{equation}
where
\begin{Eqnarray}
\label{eq:Th4f}
\theta_0 & = & - \ii \momO+\ii  \Lf{\psi}{1}{1}+\Frac{1}{4} \ii \dx^4 \mom{2}{1}, \\
\nonumber \theta_1 & = & -  2 \dx \mom{1}{1} +\Frac{1}{6} \Lf{\varphi_1+\phi_1}{1}{2} + \Frac{1}{6} \Lf{\varphi_2+\phi_2}{2}{1} ,\\
\nonumber \theta_2 & = & \ii \dt + 2\ii \dx^2 \mom{2}{1} ,\\
\nonumber \theta_3 & = & \Frac{4}{3}  \dx^3 \mom{3}{1}.
\end{Eqnarray}

Written in this form, it is clearly evident that $\Theta_4$ is free of nested
commutators and is composed of a very small number of anti-commutators. In
fact, the number of anti-commutators grows linearly with the order of
accuracy. As we see in the following section, this makes a significant
difference to the cost of our methods when compared to standard
Magnus--Lanczos schemes (which feature nested commutators and consequently a
cost that grows exponentially with order).

\begin{remark} Note that, for potentials of the form $V(x,t) =
V_0(x) + f(t)x$, the terms involving $\mom{1}{1},\mom{2}{1}$ and $\mom{3}{1}$
all vanish.
\end{remark}


\section{Implementation}

In the previous section we proposed the Magnus expansion with simplified commutators. The
next step consists in numerically approximating the exponential of this
expansion (\ref{eq:mag_disc}), which is a challenging problem itself. We will
show how Lanczos iterations can be much cheaper when combined with the
proposed versions of the Magnus expansion. In this section we present
some details of implementation and highlight some crucial features of our
schemes.


\subsection{Spatial discretisation.}
\label{sec:disc} In principle, our methods can be combined with any spatial
discretisation strategy, provided the discretisation of $\dx^n$ is symmetric
for even $n$ and skew-symmetric for odd $n$. Here we resort to spectral
collocation due to its high accuracy. Having imposed periodic boundaries on
$I$, we use equispaced grids with $N$ points. Since we work with values at
the grid points, multiplication by the function $V$ (or, in general, function
$f$) is discretised as an $N \times N$ diagonal matrix $\DDD_V$ (or $\DDD_f$)
with values of $V$ (or $f$) at the grid points along the diagonal.
The differentiation matrices $\K_k$ are symmetric for even $k$ and
skew-symmetric for odd $k$, just as we have assumed throughout. Additionally,
spectral collocation results in $\K_k$ being an $N \times N$ circulant.
Consequently, it is diagonalisable via Fourier transform,
\[ \K_k = \FFF^{-1} \DDD_{c_k} \FFF, \]
where $c_k$ is the symbol of $\K_k$ and $\FFF$ is the $N \times N$ Fourier transform matrix.

Since $ \| \DDD_f \|_2 \leq \|f\|_{\infty}$, the matrix $\DDD_f$ does not
scale with $N$. On the other hand, it can be verified that $\K_k$ scales as
$N^k$. As previously noted in \R{eq:angbradisc}, the operator $\ang{f}{k}$ is
discretised as $\Frac12 \left( \DDD_f \K_k + \K_k \DDD_f \right)$.
Consequently, upon discretisation, $\ang{f}{k}$ also scales as $N^k$. We
write $\K_k = \O{N^k}$ and, abusing notation somewhat, $\ang{f}{k} =
\O{N^k}$.

The order-four Magnus expansion $\Theta_2(h)$, with a local error $\O{\dt^5}$, discretises to the form
\begin{equation}
\label{eq:Th2disc}
\Theta_2(h) \leadsto  - \ii   \DDD_{\momO} - \left( \DDD_{\dx \mom{1}{1}} \K_1 + \K_1 \DDD_{\dx \mom{1}{1}}\right) + \ii \dt  \K_2,
\end{equation}
while the discretisation of the order-six Magnus expansion $\Theta_4(h)$ with
local error $\O{\dt^7}$ is
\begin{equation}
\label{eq:Th4disc}
\Theta_4(h) \leadsto \DDD_{\theta_0} + \Frac12 \left( \DDD_{\theta_1} \K_1 + \K_1 \DDD_{\theta_1} \right) + \Frac12 \left( \DDD_{\theta_2} \K_2 + \K_2 \DDD_{\theta_2} \right) + \Frac12 \left( \DDD_{\theta_3} \K_3 + \K_3 \DDD_{\theta_3} \right).
\end{equation}

\subsection{Evaluation of integrals and derivatives of the potential.}
\label{sec:quadrature}
Before we implement \R{eq:Th2disc} or \R{eq:Th4disc} in a practical algorithm, however, we are still left with the task of approximating functions such as $\momO$,\ $\dx \mom{1}{1}$,\ $\Lf{\psi}{1}{1}$ and $\Lf{\varphi_1+\phi_1}{1}{2}$ at the grid points, which are hidden in $\theta_i$ in the case of $\Theta_4$. These feature both integrals and derivatives of the potential. 
In some cases, it might be possible to evaluate some or all of these analytically. In other cases, however, these can be approximated by a combination of quadrature methods and finite difference differentiation\footnote{Here we suggest finite differences instead of spectral collocation since the potential is usually less oscillatory and more easily resolved than the wave function. Moreover, as we see shortly, lower degrees of accuracy are required in some cases, allowing us to reduce costs.}.

We note that since the derivatives and the integrals are in space and time, respectively, they can be exchanged. Thus the optimal strategy might involve evaluating derivatives first in some cases and integrals first in others. The optimal strategy could also depend on the relative resolutions of temporal and spatial grids.
A more challenging scenario is when the temporal grid is coarser than the spatial grid, $h = (\Delta x)^\sigma$ for some $0 < \sigma \leq 1$ (in other words, we consider larger time steps).
For the sake of simplicity, we follow a fixed strategy of evaluating the derivatives of the potential first.

{\bf Derivatives.} The various derivatives of $V$ that we need to approximate
here, $\dx V$, $\dx^2 V$, $\dx^3 V$ and $\dx^4 V$, require differentiation to
different degrees of accuracy -- by tailoring this accuracy to the $\O{h^7}$
accuracy of our Magnus expansion we can achieve the required accuracy at a
low cost.

Consider $\dx \mom{1}{1}$. Due to \R{eq:IntB},
$\Int{\zeta}{0}{\dt}{\tilde{B}_{1}(\dt,\zeta) f(\zeta) }$ is $\O{\dt^3}$ for
any $f$. Let $\K^{\mathrm{FD},n}_{k}$ be the finite difference
differentiation matrix approximating $\dx^k$ up to an error of $(\Delta
x)^n$. This accuracy, expressed in powers of $h$, naturally depends on the
relative sizes of $h$ and $\Delta x$, the latter of which is assumed to be
fixed. For instance, when $\dt = \Delta x$, it suffices to approximate $\dx
V$ to an accuracy of $\O{(\Delta x)^4} = \O{\dt^4}$ via
$\K^{\mathrm{FD},4}_{1} \MM{V}$, since the integral
$\Int{\zeta}{0}{\dt}{\tilde{B}_{1}(\dt,\zeta) \K^{\mathrm{FD},4}_{1}
\MM{V}(\zeta)}$ approximates $\dx \mom{1}{1}$ to the required accuracy of
$\O{\dt^7}$. When the time step is larger, say $\dt = \sqrt{\Delta x}$, the
lower accuracy (and lower cost) differentiation matrix
$\K^{\mathrm{FD},2}_{1}$ suffices. A practical method could use
$\K^{\mathrm{FD},4}_{1}$ for $\Delta x \leq h < \sqrt{\Delta x}$ and
$\K^{\mathrm{FD},2}_{1}$ for $\sqrt{\Delta x} \leq h$.

Similar considerations show that we need to approximate $\dx^2 V$ to an accuracy of
$\O{h^3}$, $\dx^3 V$ to an accuracy of $\O{h^2}$ and $\dx^4 V$ to an accuracy of $\O{h^3}$.

{\bf Quadrature.} For the purpose of approximating the integrals, we can
resort to a variety of quadrature methods, among which Gauss--Legendre
quadratures are the most popular due to their high orders of accuracy. For
instance, all these integrals can be approximated to $\O{\dt^7}$ accuracy
using Gauss--Legendre quadrature at the knots $\tau_k = h (1 + k
\sqrt{3/5})/2 ,\ k=-1,0,1$, with the weights $w_k = \Frac{5}{18}\dt, \Frac49
\dt, \Frac{5}{18} \dt$ \cite{davis84quad}. \footnote{Recall that since these
Magnus expansions are odd in $h$, the  $\O{\dt^6}$ Gauss-Legendre quadrature
automatically becomes $\O{\dt^7}$ in this context.}


%
%
%

{\bf Approximation of line integrals.} Under $\sigma = 1$, the line integrals $\mom{j}{k}$ and their derivatives appearing in $\Theta_4$ can be approximated to $\O{\dt^7}$ accuracy by using the weights $w_k$,
\begin{Eqnarray}
\label{eq:FDu1}
    \momO & \leadsto & w_{-1} \MM{V}(\tau_{-1}) + w_{0} \MM{V}(\tau_{0})  + w_{1} \MM{V}(\tau_{1}),\\
\label{eq:FDu2}
    \dx \mom{1}{1} & \leadsto & w_{-1} \tilde{B}_{1}(\dt,\tau_{-1}) \K^{\mathrm{FD},4}_{1} \MM{V}(\tau_{-1})  +  w_{1} \tilde{B}_{1}(\dt,\tau_{1}) \K^{\mathrm{FD},4}_{1} \MM{V}(\tau_{1}),\\
    \dx^2 \mom{2}{1} & \leadsto & w_{-1} \tilde{B}_{2}(\dt,\tau_{-1}) \K^{\mathrm{FD},3}_{2} \MM{V}(\tau_{-1})  \nonumber  \\
\label{eq:FDu3}
    & & +  w_{0} \tilde{B}_{2}(\dt,\tau_{0}) \K^{\mathrm{FD},3}_{2} \MM{V}(\tau_{0}) + w_{1} \tilde{B}_{2}(\dt,\tau_{1}) \K^{\mathrm{FD},3}_{2} \MM{V}(\tau_{1}),\\
    \dx^3 \mom{1}{3} & \leadsto & w_{-1} \tilde{B}_{1}(\dt,\tau_{-1})^3 \K^{\mathrm{FD},2}_{3} \MM{V}(\tau_{-1})  \nonumber  \\
\label{eq:FDu4}
    & & +  w_{0} \tilde{B}_{1}(\dt,\tau_{0})^3 \K^{\mathrm{FD},2}_{3} \MM{V}(\tau_{0}) + w_{1} \tilde{B}_{1}(\dt,\tau_{1})^3 \K^{\mathrm{FD},2}_{3} \MM{V}(\tau_{1}),\\
\label{eq:FDu5}
    \dx^4 \mom{2}{1} & \leadsto & w_{-1} \tilde{B}_{2}(\dt,\tau_{-1}) \K^{\mathrm{FD},3}_{4} \MM{V}(\tau_{-1})  \nonumber  \\
    & & +  w_{0} \tilde{B}_{2}(\dt,\tau_{0}) \K^{\mathrm{FD},3}_{4} \MM{V}(\tau_{0}) + w_{1} \tilde{B}_{2}(\dt,\tau_{1}) \K^{\mathrm{FD},3}_{4} \MM{V}(\tau_{1}),
\end{Eqnarray}
where we note that since $\tilde{B}_{1}(\dt,\tau_{0}) =0$, the $\tau_0$ term
does not appear in \R{eq:FDu2}. We note that, instead of \R{eq:FDu3} and
\R{eq:FDu5}, approximating $\mom{2}{1}$ first and then evaluating its
derivatives would be cheaper overall. However, as mentioned earlier,  we
attempt here to provide a simple and clear procedure, not a fully optimised
one.

For the order four Magnus expansion, $\Theta_2$, the first two terms  \R{eq:FDu1} and \R{eq:FDu2} suffice. However, since we need only $\O{h^5}$ accuracy, we could do with just two Gauss--Legendre knots.

{\bf Approximation of integrals over the triangle.} For the integrals over the triangle such as $\Lf{\psi}{1}{1}$, the appropriate weights can be found by substituting the interpolant, $\tilde{\MM{v}}(t) =  \sum_{k=-1}^1 \ell_{k}(t) \MM{v}(\tau_k)$, where $\MM{v}$ is usually a derivative of the potential\footnote{For instance, $\MM{v}(t) = \K^{\mathrm{FD},3}_{1} \MM{V}(t)$ suffices in the approximation of $\Lf{\psi}{1}{1}$.} and where $\ell_{k}(t)$ are the Lagrange cardinal functions, $\ell_{k}(\tau_j) = \delta_{j,k}$. Thus we discretise,
\begin{Eqnarray}
\nonumber \Lf{f}{a}{b} &\leadsto& \Int{\zeta}{0}{\dt}{\Int{\xi}{0}{\zeta}{\sum_{j=-1}^1 \sum_{k=-1}^1 f(\dt,\zeta,\xi) \ell_j(\zeta) \ell_k(\xi) \left[ \K^{\mathrm{FD},r}_{a} \MM{V}(\tau_j)\right] \left[\K^{\mathrm{FD},r}_{b} \MM{V}(\tau_k) \right] }}\\
\label{eq:Lf} &=& \sum_{j=-1}^1 \sum_{k=-1}^1 w_{jk}^f \left[
\K^{\mathrm{FD},r}_{a} \MM{V}(\tau_j)\right] \left[\K^{\mathrm{FD},r}_{b}
\MM{V}(\tau_k) \right]\!,
\end{Eqnarray}
where we need to approximate derivatives of $V$ to order $r$ (under the
scaling $\sigma = 1$, $r=3$ suffices for all $\Lambda$ terms in $\Theta_4$
for $\O{h^7}$ accuracy), and where $w_{jk}^f$ are the weights specific to
$f$,
\[w_{jk}^f = \Int{\zeta}{0}{\dt}{\Int{\xi}{0}{\zeta}{ f(\dt,\zeta,\xi) \ell_j(\zeta) \ell_k(\xi) }}. \]
The weights for the functions $\psi,\varphi_1,\varphi_2,\phi_1$ and $\phi_2$
that are required for the implementation of an order-six method have been
provided in
 Appendix~\ref{sec:quadratureweights}. The method for
discretising the anti-commutators, as well as a particular recipe for
approximating the integrals and derivatives of the potential, is in place.

\begin{remark}Having elaborated on the use of Gauss--Legendre quadratures in developing a specific scheme,
we remind the reader that a major advantage of
preserving integrals throughout the workings in Section~2 is the flexibility
of allowing alternative means for evaluating integrals and derivatives,
including the possibility of exact integration and derivation.
\end{remark}

\subsection{Approximation of the exponential of a Magnus expansion.}
\label{sec:exponential} After discretising $\Theta_2(h)$ and $\Theta_4(h)$,
we are left with the task of approximating their exponential in
\R{eq:magnus_exp_t} in order to find the solution
\[ \MM{u}^1 = \exp(\Theta_m(h)) \MM{u}^0. \]
As discussed in Subsection~\ref{sec:existing}, Lanczos iterations are a very
effective, and perhaps the most popular, means for approximating the
exponential of a Magnus expansion. This will be the approach adopted in this
paper.

Approximation of the matrix vector product $\exp(\Theta_p) \MM{u}$ via Lanczos iterations requires the evaluation of $\Theta_p \MM{v}$ in each Lanczos iteration. So long as the number of steps is reasonably small, the cost is dominated by the cost of evaluating $\Theta_p \MM{v}$.

Standard Magnus expansions feature nested commutators. In the method
presented in Subsection~\ref{sec:existing}, $\Theta_p$ features commutators
nested to $p-1$ levels. For a commutator $C_p$ that is nested to $p$ levels,
the cost of evaluating the matrix-vector product $C_p \MM{v}$ grows
exponentially with $p$. Consequently, the cost of $\Theta_p \MM{v}$ in
standard Magnus--Lanczos schemes grows exponentially with the order of the
scheme. In contrast our proposed Magnus expansions feature a linearly growing
number of non-nested terms. As evident from \R{eq:Th4disc}, for instance our
$\O{h^7}$ method $\Theta_4(h)$ is comprised of a small number of terms. The
approximation of $\Theta_4 \MM{v}$ in each Lanczos iterations requires us to
compute
 \begin{Eqnarray}
\label{eq:Th4v}
\nonumber \Theta_4(h)\MM{v} & = & \DDD_{\theta_0} \MM{v} + \Frac12 \left( \DDD_{\theta_1} \K_1 + \K_1 \DDD_{\theta_1} \right) \MM{v} + \Frac12 \left( \DDD_{\theta_2} \K_2 + \K_2 \DDD_{\theta_2} \right) \MM{v} + \Frac12 \left( \DDD_{\theta_3} \K_3 + \K_3 \DDD_{\theta_3} \right)\MM{v}\\
\nonumber & = & \DDD_{\theta_0} \MM{v} + \Frac12 \left( \DDD_{\theta_1} \FFF^{-1} \DDD_{c_1} \FFF + \FFF^{-1} \DDD_{c_1} \FFF \DDD_{\theta_1} \right)\MM{v} \\
\nonumber & & + \Frac12 \left( \DDD_{\theta_2} \FFF^{-1} \DDD_{c_2} \FFF + \FFF^{-1} \DDD_{c_2} \FFF \DDD_{\theta_2} \right) \MM{v} + \Frac12 \left( \DDD_{\theta_3} \FFF^{-1} \DDD_{c_3} \FFF + \FFF^{-1} \DDD_{c_3} \FFF \DDD_{\theta_3} \right)\MM{v}\\
& = & \DDD_{\theta_0} \MM{v} + \Frac12 \left( \sum_{j=1}^3 \DDD_{\theta_j} \FFF^{-1} \DDD_{c_j} \right) \FFF \MM{v} + \Frac12 \FFF^{-1} \left( \sum_{j=1}^3   \DDD_{c_j} \FFF \DDD_{\theta_j} \MM{v}\right),
\end{Eqnarray}
which requires merely eight FFTs.

\subsection{Unitarity, norm preservation and stability}
Note that our Magnus expansions \R{eq:Omt3} and \R{eq:Th4} are of the form $
\sum_{k=0}^\infty \ii^{k+1} \Ang{k}{\theta_k} $ for some $\theta_k$, and it
can be seen that each term $\ii^{k+1} \Ang{k}{\theta_k}$ discretises to a
skew-Hermitian form in \R{eq:Th2disc} and \R{eq:Th4disc}, respectively. The
Magnus expansion, developed in this way preserves skew-Hermiticity and its
exponential therefore preserves unitarity. As mentioned in
Subsection~\ref{sec:existing}, this is consistent with a central feature of
quantum mechanics. Additionally, since the exponential is unitary,
\[ \| \MM{u}^{1} \|_{2} = \| \exp\left(\Theta_m(\dt)\right) \MM{u}^0\|_{2} = \| \MM{u}^0\|_{2}, \]
and the norm of $\MM{u}$ is preserved. Consequently, unitarity guarantees
stability of our schemes under any scaling of $h$ and $\Delta x$.

\section{Numerical Examples}
\label{sec:experiments}

\begin{figure}[tbh]
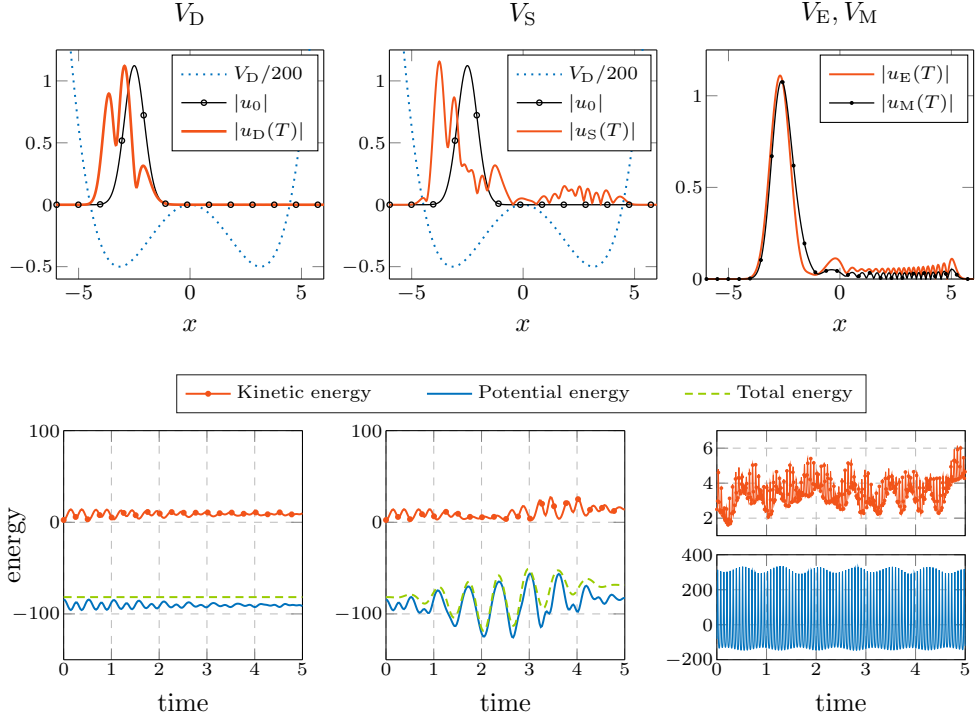

        \input{uDWInitialFinal.tex}
        \input{uS10InitialFinal.tex}
        \input{uE100InitialFinal.tex}
        \\
        \\
%
\definecolor{mycolor1}{rgb}{0.00000,0.44706,0.74118}%
\definecolor{mycolor2}{rgb}{0.95000,0.32500,0.09800}%
\definecolor{mycolor3}{rgb}{1.00000,0.99000,0.97000}%
\begin{tikzpicture}
    \begin{axis}[%
    width=5.25in, height=0.75in,
    hide axis,
    xmin=10,
    xmax=50,
    ymin=0,
    ymax=0.4,
    legend style={draw=white!15!black,legend cell align=left, legend columns=-1, font=\scriptsize}
    ]
    \addlegendimage{mycolor2, line width=0.75pt, mark=*, mark size=0.75, mark options={solid}}
    \addlegendentry{Kinetic energy $\qquad$};
    \addlegendimage{mycolor1, line width=0.75pt}
    \addlegendentry{Potential energy $\qquad$};
    \addlegendimage{black!20!lime, line width=0.75pt, densely dashed}
    \addlegendentry{Total energy};
    \end{axis}
\end{tikzpicture}\\
        \input{energyDW.tex}
        \input{energyS10.tex}
        \input{energyE100.tex}
    \caption{{\bf [top row]} The initial condition $u_0$ evolves to $u_{\mathrm{D}}$ under the influence of $V_{\mathrm{D}}$ (left),
    to $u_{\mathrm{S}}$ under $V_{\mathrm{S}}$ (centre) and to $u_{\mathrm{E}}$ and $u_{\mathrm{M}}$ under $V_{\mathrm{E}}$
    and $V_{\mathrm{M}}$, respectively (right). The potential $V_{\mathrm{D}}$ is scaled down for ease of presentation. {\bf [bottom row]} Corresponding evolution of energies.}
    \label{fig:uIF}
\end{figure}

The initial condition for our numerical experiments is a Gaussian wavepacket
\[u_0(x) = (\delta  \pi)^{-1/4}  \exp\left( (-(x - x_0)^2 )/ (2 \delta)\right),\quad x_0 = -2.5,\quad \delta =0.2,\]
sitting in the left well of a double well potential,
\[V_{\mathrm{D}}(x) = x^4 - 20 x^2. \]
We take $[-10,10]$ as our spatial domain and $[0,5]$ as our temporal domain.
When we allow the wave function to evolve under $V_{\mathrm{D}}$, it remains
largely confined to the left well at the final time, $T=5$ (see
Figure~\ref{fig:uIF}, top left).
Superimposing a time dependent excitation to the potential, we are able to
exert control on the wave function. In  Figure~\ref{fig:uIF} (top centre and
right) we show the influence of two excitations of the form $f(t)\,x$ with
the choices $f(t) = 10\, \mathrm{S}_{10,T}(t)$ and $f(t) = -25\,
\mathrm{E}_{100}(t)$, where
\[ \mathrm{S}_{\omega,T}(t) = \sin((\pi t/T)^2) \sin(\omega t), \quad \mathrm{E}_{\omega,T}(t) = \exp(2 \sin(\omega t))-1.\]
The effective  time-varying potentials in these cases are
\[V_{\mathrm{S}}(x,t) = V_{\mathrm{D}}(x) + 10\, \mathrm{S}_{10,T}(t) x,\quad V_{\mathrm{E}}(x,t) = V_{\mathrm{D}}(x) - 25\, \mathrm{E}_{100}(t) x,\]
respectively. Since the potentials are available in their analytic form, we
use analytic derivatives in our implementation. The integrals were
approximated via three Gauss--Legendre knots, as outlined in
Subsection~\ref{sec:quadrature}. In principle we can also use analytic or
asymptotic approximations for the integrals. Potential accuracy advantages of
resorting to analytic approximations should become evident by comparing with
a higher degree quadrature -- for this purpose we also present results using
eleven Gauss--Legendre knots.

{\bf Mean field approximation ($V_{\mathrm{M}}$)}. The excitation in
$V_{\mathrm{S}}$ features a low frequency oscillation at $\omega = 10$, while
$V_{\mathrm{E}}$ has a higher frequency oscillation, $\omega = 100$. In the
limit $\omega \rightarrow \infty$ the effect of the potential function can be
approximated by a mean field potential\footnote{Effectively the first and
trivial Magnus expansion $\Theta_1$.} and it is worth finding out to what
extent this approximation suffices for $V_{\mathrm{E}}$. Since
$\Int{t}{0}{5}{\mathrm{E}_{100}(t)} \approx 6.45083$, the (time-independent)
mean field potential is roughly
\[V_{\mathrm{M}}(x) = V_{\mathrm{DW}}(x) -32.25415 x.\] In Figure~\ref{fig:uIF} (top right) it is evident that the mean field
solution $u_{\mathrm{M}}(T)$ is not a sufficiently accurate approximation to
$u_{\mathrm{E}}(T)$, and at $\omega=100$ we do require a solution via
high-order Magnus based methods.


{\bf Magnus--\texttt{expm}}. In the numerical experiments presented in this
section, order four and six traditional Magnus expansions are denoted by
$\mathrm{M_4}$ and $\mathrm{M_6}$, respectively, while the corresponding
 Magnus expansions with simplified commutators are denoted by $\mathrm{S_4}$ and
$\mathrm{S_6}$, respectively. The order four methods use two Gauss--Legendre
quadrature knots, while order six methods use three knots. All these methods
use $180$ spatial grid points and are exponentiated using \textsc{Matlab}'s
\texttt{expm}. We present the errors for these Magnus--\texttt{expm} methods
in order to study the error inherent in the Magnus expansion separately from
the error due to Lanczos iterations.

{\bf Higher accuracy quadratures}. In a high frequency regime, we encounter
more oscillations per time step and three quadrature knots can be inadequate
for approximating the integrals adequately. In such cases, we can expect to
see a considerable advantage when using analytic integrals, asymptotic
approximations or higher accuracy quadratures. This behaviour is exhibited
in Figure~\ref{fig:S6M6VE}, where our order six integral preserving Magnus
expansion $\mathrm{S_6}$ is seen to have a much higher accuracy than the
standard Magnus expansion $\mathrm{M_6}$, particularly when combined with a
higher accuracy approximation to the integrals. In this case, we resort to 11
Gauss--Legendre quadratures, denoted by the postfix $\mathrm{G_{11}}$.
Analytic integrals could possibly improve the accuracy further, as could
highly oscillatory quadrature methods \cite{DHIbook}, which can be easily
transplanted to this setting.

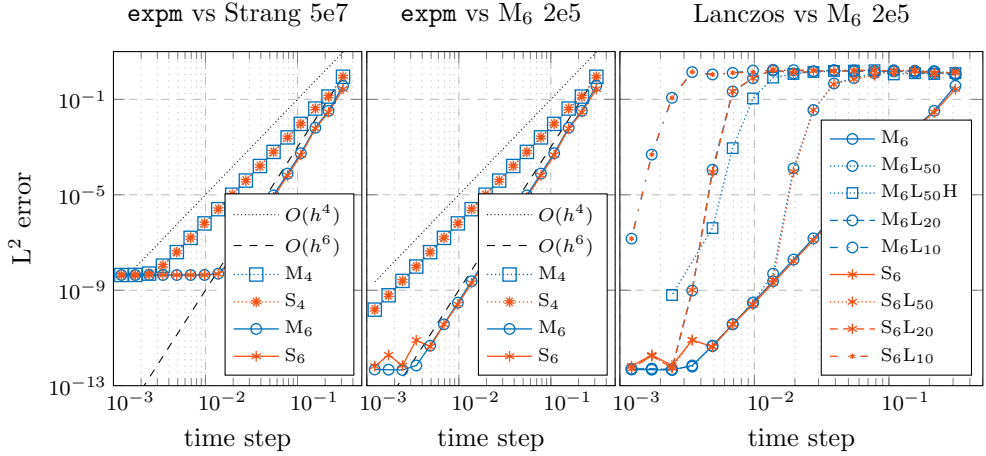
\begin{figure}[tbh]
%
%
\definecolor{mycolor1}{rgb}{0.00000,0.44700,0.74100}%
\definecolor{mycolor2}{rgb}{0.85000,0.32500,0.09800}%
\definecolor{mycolor3}{rgb}{0.92900,0.69400,0.12500}%
\definecolor{mycolor4}{rgb}{0.49400,0.18400,0.55600}%
\begin{tikzpicture}

\begin{axis}[%
width=1.3in, height=1.75in, scale only axis, xmode=log, xmin=0.001, xmax=0.5,
xminorticks=true, ymode=log, ymin=1e-13, ymax=10, yminorticks=true, axis
background/.style={fill=white}, xticklabel style = {xshift=1.5ex},
xlabel={time step},  ylabel={$\mathrm{L}^2$ error}, title={\texttt{expm} vs
Strang 5e7},legend style={legend cell align=left, align=left,
draw=black,font=\scriptsize}, legend pos=south east, grid=both]

\addplot [color=black, densely dotted]
  table[row sep=crcr]{%
0.3125	9.5367431640625\\
0.217391304347826	2.23341111559779\\
0.15625	0.596046447753906\\
0.108695652173913	0.139588194724862\\
0.078125	0.0372529029846191\\
0.0549450549450549	0.00911411382355726\\
0.0390625	0.0023283064365387\\
0.0274725274725275	0.000569632113972329\\
0.01953125	0.000145519152283669\\
0.0137741046831956	3.5995939830998e-05\\
0.009765625	9.09494701772928e-06\\
0.00689655172413793	2.26218433691842e-06\\
0.0048828125	5.6843418860808e-07\\
0.00345065562456867	1.41777226236882e-07\\
0.00244140625	3.5527136788005e-08\\
0.00172592336900242	8.87331780757525e-09\\
0.001220703125	2.22044604925031e-09\\
}; \addlegendentry{$O(h^4)$}

\addplot [color=black, dashed]
  table[row sep=crcr]{%
0.3125	0.931322574615479\\
0.217391304347826	0.105548729470595\\
0.15625	0.0145519152283669\\
0.108695652173913	0.00164919889797804\\
0.078125	0.000227373675443232\\
0.0549450549450549	2.7515136528068e-05\\
0.0390625	3.5527136788005e-06\\
0.0274725274725275	4.29924008251063e-07\\
0.01953125	5.55111512312578e-08\\
0.0137741046831956	6.82936423418975e-09\\
0.009765625	8.67361737988404e-10\\
0.00689655172413793	1.07594974407535e-10\\
0.0048828125	1.35525271560688e-11\\
0.00345065562456867	1.68814486939283e-12\\
0.00244140625	2.11758236813575e-13\\
0.00172592336900242	2.64319409124601e-14\\
0.001220703125	3.30872245021211e-15\\
};\addlegendentry{$O(h^6)$}

\addplot [color=colorclassyblue,line width=0.5pt, densely dotted,
mark=square,mark size=2.5,mark options={solid},mark repeat=1]
  table[row sep=crcr]{%
0.3125	0.882928921771494\\
0.217391304347826	0.130522094333576\\
0.15625	0.0411017276105357\\
0.108695652173913	0.00936691054235059\\
0.078125	0.00254447988558078\\
0.0549450549450549	0.000627631793894988\\
0.0390625	0.000160943822128434\\
0.0274725274725275	3.94513716247398e-05\\
0.01953125	1.00875738382437e-05\\
0.0137741046831956	2.496443581155e-06\\
0.009765625	6.30905210688566e-07\\
0.00689655172413793	1.56984324403961e-07\\
0.0048828125	3.96511630233368e-08\\
0.00345065562456867	1.0725186443457e-08\\
0.00244140625	4.98570608752184e-09\\
0.00172592336900242	4.38892038517723e-09\\
0.001220703125	4.35096462269466e-09\\
}; \addlegendentry{$\mathrm{M_4}$}

\addplot [color=colorclassyorange, line width=0.5pt, densely dotted,
mark=10-pointed star,mark size=2,mark options={solid},mark repeat=1]
  table[row sep=crcr]{%
0.3125	0.882928921771639\\
0.217391304347826	0.130522094333446\\
0.15625	0.0411017276105277\\
0.108695652173913	0.00936691054230564\\
0.078125	0.00254447988557402\\
0.0549450549450549	0.00062763179386562\\
0.0390625	0.000160943822142628\\
0.0274725274725275	3.9451371557971e-05\\
0.01953125	1.00875738788327e-05\\
0.0137741046831956	2.496443521078e-06\\
0.009765625	6.30905241075064e-07\\
0.00689655172413793	1.56984295756559e-07\\
0.0048828125	3.96511662554658e-08\\
0.00345065562456867	1.0725193395233e-08\\
0.00244140625	4.98573744273586e-09\\
0.00172592336900242	4.3889092818853e-09\\
0.001220703125	4.35096264830559e-09\\
}; \addlegendentry{$\mathrm{S_4}$}

\addplot [color=colorclassyblue,line width=0.5pt,mark=o,mark size=2,mark
options={solid},mark repeat=1]
  table[row sep=crcr]{%
0.3125	0.368201464899316\\
0.217391304347826	0.0329875422536024\\
0.15625	0.00644972320381661\\
0.108695652173913	0.000552444373393711\\
0.078125	7.73343900108409e-05\\
0.0549450549450549	9.48477358224664e-06\\
0.0390625	1.23139344781019e-06\\
0.0274725274725275	1.49488523565676e-07\\
0.01953125	1.97912594262925e-08\\
0.0137741046831956	4.95020529890938e-09\\
0.009765625	4.35869856654883e-09\\
0.00689655172413793	4.3492445178785e-09\\
0.0048828125	4.34919094049515e-09\\
0.00345065562456867	4.34920247628818e-09\\
0.00244140625	4.34921267124911e-09\\
0.00172592336900242	4.34919445428285e-09\\
0.001220703125	4.34919381050291e-09\\
}; \addlegendentry{$\mathrm{M_6}$}

\addplot [color=colorclassyorange, line width=0.5pt,mark=asterisk,mark
size=2,mark options={solid},mark repeat=1]
  table[row sep=crcr]{%
0.3125	0.26984508613968\\
0.217391304347826	0.0287740902917036\\
0.15625	0.00598535901910583\\
0.108695652173913	0.000464899336300946\\
0.078125	6.49508497535545e-05\\
0.0549450549450549	7.97745829554796e-06\\
0.0390625	1.03631057786135e-06\\
0.0274725274725275	1.25862502017352e-07\\
0.01953125	1.6822917279445e-08\\
0.0137741046831956	4.78156889213143e-09\\
0.009765625	4.35561539704268e-09\\
0.00689655172413793	4.34923854697655e-09\\
0.0048828125	4.34919205687439e-09\\
0.00345065562456867	4.34924675343942e-09\\
0.00244140625	4.34922300482269e-09\\
0.00172592336900242	4.3492202783684e-09\\
0.001220703125	4.34919158962728e-09\\
}; \addlegendentry{$\mathrm{S_6}$}

\end{axis}
\end{tikzpicture}%
        \hspace{-0.35cm}
%
%
\definecolor{mycolor1}{rgb}{0.00000,0.44700,0.74100}%
\definecolor{mycolor2}{rgb}{0.85000,0.32500,0.09800}%
\definecolor{mycolor3}{rgb}{0.92900,0.69400,0.12500}%
\definecolor{mycolor4}{rgb}{0.49400,0.18400,0.55600}%
\definecolor{mycolor5}{rgb}{0.46600,0.67400,0.18800}%
\definecolor{mycolor6}{rgb}{0.30100,0.74500,0.93300}%
\begin{tikzpicture}

\begin{axis}[%
width=1.3in, height=1.75in, scale only axis, xmode=log, xmin=0.001, xmax=0.5,
xminorticks=true, ymode=log, ymin=1e-13, ymax=10, yminorticks=true, axis
background/.style={fill=white}, xticklabel style = {xshift=1.5ex},
yticklabels={},xlabel={time step},title={\texttt{expm} vs $\mathrm{M_6}$
2e5},legend style={legend cell align=left, align=left,
draw=black,font=\scriptsize}, legend pos=south east, grid=both]

\addplot [color=black, densely dotted]
  table[row sep=crcr]{%
  0.3125	9.5367431640625\\
0.217391304347826	2.23341111559779\\
0.15625	0.596046447753906\\
0.108695652173913	0.139588194724862\\
0.078125	0.0372529029846191\\
0.0549450549450549	0.00911411382355726\\
0.0390625	0.0023283064365387\\
0.0274725274725275	0.000569632113972329\\
0.01953125	0.000145519152283669\\
0.0137741046831956	3.5995939830998e-05\\
0.009765625	9.09494701772928e-06\\
0.00689655172413793	2.26218433691842e-06\\
0.0048828125	5.6843418860808e-07\\
0.00345065562456867	1.41777226236882e-07\\
0.00244140625	3.5527136788005e-08\\
0.00172592336900242	8.87331780757525e-09\\
0.001220703125	2.22044604925031e-09\\
}; \addlegendentry{$O(h^4)$}

\addplot [color=black, dashed]
  table[row sep=crcr]{%
0.3125	0.931322574615479\\
0.217391304347826	0.105548729470595\\
0.15625	0.0145519152283669\\
0.108695652173913	0.00164919889797804\\
0.078125	0.000227373675443232\\
0.0549450549450549	2.7515136528068e-05\\
0.0390625	3.5527136788005e-06\\
0.0274725274725275	4.29924008251063e-07\\
0.01953125	5.55111512312578e-08\\
0.0137741046831956	6.82936423418975e-09\\
0.009765625	8.67361737988404e-10\\
0.00689655172413793	1.07594974407535e-10\\
0.0048828125	1.35525271560688e-11\\
0.00345065562456867	1.68814486939283e-12\\
0.00244140625	2.11758236813575e-13\\
0.00172592336900242	2.64319409124601e-14\\
0.001220703125	3.30872245021211e-15\\
}; \addlegendentry{$O(h^6)$}

\addplot [color=colorclassyblue,line width=0.5pt, densely dotted,
mark=square,mark size=2.5,mark options={solid},mark repeat=1]
  table[row sep=crcr]{%
0.3125	0.882928919849092\\
0.217391304347826	0.130522094072298\\
0.15625	0.0411017275475281\\
0.108695652173913	0.00936691054876091\\
0.078125	0.00254447990729026\\
0.0549450549450549	0.000627631819860998\\
0.0390625	0.0001609438490817\\
0.0274725274725275	3.9451398663124e-05\\
0.01953125	1.00876002420156e-05\\
0.0137741046831956	2.496467149862e-06\\
0.009765625	6.30917580923116e-07\\
0.00689655172413793	1.56951428489437e-07\\
0.0048828125	3.94392869554529e-08\\
0.00345065562456867	9.83118549441378e-09\\
0.00244140625	2.46508868227278e-09\\
0.00172592336900242	6.16733943084453e-10\\
0.001220703125	1.54106623368269e-10\\
}; \addlegendentry{$\mathrm{M_4}$}

\addplot [color=colorclassyorange, line width=0.5pt, densely dotted,
mark=10-pointed star,mark size=2,mark options={solid},mark repeat=1]
  table[row sep=crcr]{%
0.3125	0.882928919849237\\
0.217391304347826	0.130522094072168\\
0.15625	0.04110172754752\\
0.108695652173913	0.00936691054871597\\
0.078125	0.00254447990728348\\
0.0549450549450549	0.000627631819831638\\
0.0390625	0.000160943849095902\\
0.0274725274725275	3.94513985963519e-05\\
0.01953125	1.00876002825777e-05\\
0.0137741046831956	2.49646708977115e-06\\
0.009765625	6.30917611156182e-07\\
0.00689655172413793	1.56951399745496e-07\\
0.0048828125	3.94392908465315e-08\\
0.00345065562456867	9.83118454999105e-09\\
0.00244140625	2.46509653925427e-09\\
0.00172592336900242	6.16735450248007e-10\\
0.001220703125	1.54095029165826e-10\\
}; \addlegendentry{$\mathrm{S_4}$}

\addplot [color=colorclassyblue,line width=0.5pt,mark=o,mark size=2,mark
options={solid},mark repeat=1]
  table[row sep=crcr]{%
0.3125	0.368201464108403\\
0.217391304347826	0.0329875422034794\\
0.15625	0.00644972321069644\\
0.108695652173913	0.000552444384023556\\
0.078125	7.73344040650404e-05\\
0.0549450549450549	9.48478726104196e-06\\
0.0390625	1.23140042165292e-06\\
0.0274725274725275	1.49439881005864e-07\\
0.01953125	1.93220958353901e-08\\
0.0137741046831956	2.37882719422355e-09\\
0.009765625	3.02221731424032e-10\\
0.00689655172413793	3.74751850531667e-11\\
0.0048828125	4.75110981531588e-12\\
0.00345065562456867	7.03424463236684e-13\\
0.00244140625	4.54984946705714e-13\\
0.00172592336900242	4.67771859792481e-13\\
0.001220703125	4.77919368977232e-13\\
}; \addlegendentry{$\mathrm{M_6}$}

\addplot [color=colorclassyorange, line width=0.5pt,mark=asterisk,mark
size=2,mark options={solid},mark repeat=1]
  table[row sep=crcr]{%
0.3125	0.269845085562075\\
0.217391304347826	0.0287740902541452\\
0.15625	0.00598535902887919\\
0.108695652173913	0.000464899349721709\\
0.078125	6.49508669185554e-05\\
0.0549450549450549	7.97747492861161e-06\\
0.0390625	1.03631923419894e-06\\
0.0274725274725275	1.25805095259265e-07\\
0.01953125	1.62687618157949e-08\\
0.0137741046831956	2.00462025471802e-09\\
0.009765625	2.54788425922393e-10\\
0.00689655172413793	3.72330841306388e-11\\
0.0048828125	4.31720502106715e-12\\
0.00345065562456867	8.12426190960238e-12\\
0.00244140625	7.03046271063629e-13\\
0.00172592336900242	1.94190526176589e-12\\
0.001220703125	6.61014646447738e-13\\
}; \addlegendentry{$\mathrm{S_6}$}

\end{axis}
\end{tikzpicture}%
        \hspace{-0.35cm}
%
%
\definecolor{mycolor1}{rgb}{0.00000,0.44700,0.74100}%
\definecolor{mycolor2}{rgb}{0.85000,0.32500,0.09800}%
\definecolor{mycolor3}{rgb}{0.92900,0.69400,0.12500}%
\definecolor{mycolor4}{rgb}{0.49400,0.18400,0.55600}%
\definecolor{mycolor5}{rgb}{0.46600,0.67400,0.18800}%
\definecolor{mycolor6}{rgb}{0.30100,0.74500,0.93300}%
\definecolor{mycolor7}{rgb}{0.63500,0.07800,0.18400}%
\begin{tikzpicture}

\begin{axis}[%
width=1.9in, height=1.75in, scale only axis, xmode=log, xmin=0.001, xmax=0.5,
xminorticks=true, ymode=log, ymin=1e-13, ymax=10, yminorticks=true, axis
background/.style={fill=white}, xticklabel style = {xshift=1.5ex},
xlabel={time step}, yticklabels={}, title={Lanczos vs $\mathrm{M_6}$
2e5},legend style={legend cell align=left, align=left,
draw=black,font=\scriptsize}, legend pos=south east, grid=both]

\addplot [color=colorclassyblue,line width=0.5pt,mark=o,mark size=2,mark
options={solid},mark repeat=1]
  table[row sep=crcr]{%
0.3125	0.368201464108403\\
0.217391304347826	0.0329875422034794\\
0.15625	0.00644972321069644\\
0.108695652173913	0.000552444384023556\\
0.078125	7.73344040650404e-05\\
0.0549450549450549	9.48478726104196e-06\\
0.0390625	1.23140042165292e-06\\
0.0274725274725275	1.49439881005864e-07\\
0.01953125	1.93220958353901e-08\\
0.0137741046831956	2.37882719422355e-09\\
0.009765625	3.02221731424032e-10\\
0.00689655172413793	3.74751850531667e-11\\
0.0048828125	4.75110981531588e-12\\
0.00345065562456867	7.03424463236684e-13\\
0.00244140625	4.54984946705714e-13\\
0.00172592336900242	4.67771859792481e-13\\
0.001220703125	4.77919368977232e-13\\
}; \addlegendentry{$\mathrm{M_6}$}

\addplot [color=colorclassyblue, line width=0.5pt,densely dotted,mark=o,mark
size=2,mark options={solid},mark repeat=1]
  table[row sep=crcr]{%
0.3125	1.21524034522093\\
0.217391304347826	1.53444807326375\\
0.15625	1.58731993597661\\
0.108695652173913	1.66242582192179\\
0.078125	1.2623005578685\\
0.0549450549450549	0.779029970006988\\
0.0390625	0.462685285364014\\
0.0274725274725275	0.0357017245696437\\
0.01953125	0.0001291105193747\\
0.0137741046831956	4.96439109670396e-09\\
0.009765625	3.02057392319797e-10\\
0.00689655172413793	3.73314569184889e-11\\
0.0048828125	4.56642450370452e-12\\
0.00345065562456867	6.29476944218422e-13\\
0.00244140625	4.81524133045081e-13\\
0.00172592336900242	5.13521730357047e-13\\
0.001220703125	5.04547967010372e-13\\
}; \addlegendentry{$\mathrm{M_6L_{50}}$}

\addplot [color=colorclassyblue, line width=0.5pt,densely
dotted,mark=square,mark size=2,mark options={solid},mark repeat=1]
  table[row sep=crcr]{%
0.3125	1.26665250709716\\
0.217391304347826	1.15238064949503\\
0.15625	1.2435780818431\\
0.108695652173913	1.1087211813278\\
0.078125	1.65997689194837\\
0.0549450549450549	1.60307738319895\\
0.0390625	1.5417372230802\\
0.0274725274725275	1.40760535504811\\
0.01953125	1.15535136375467\\
0.0137741046831956	0.814537536265418\\
0.009765625	0.108724505085733\\
0.00689655172413793	0.000899928278923905\\
0.0048828125	4.01931799591027e-07\\
0.00244140625	6.31896777287299e-10\\
}; \addlegendentry{$\mathrm{M_6L_{50}H}$}

\addplot [color=colorclassyblue, line width=0.5pt,densely dashed,mark=o,mark
size=2,mark options={solid},mark repeat=1]
  table[row sep=crcr]{%
0.3125	1.12300699608571\\
0.217391304347826	1.40545436453362\\
0.15625	1.42852218989758\\
0.108695652173913	1.48478686199886\\
0.078125	1.51990705489631\\
0.0549450549450549	1.47723098515769\\
0.0390625	1.60745600320055\\
0.0274725274725275	1.49405089009665\\
0.01953125	1.21676496470502\\
0.0137741046831956	1.72887755236166\\
0.009765625	0.759719376091832\\
0.00689655172413793	0.21434432851834\\
0.0048828125	0.000111836772568607\\
0.00345065562456867	9.69655186864352e-10\\
0.00244140625	4.80956675618052e-13\\
0.00172592336900242	5.12519133979299e-13\\
0.001220703125	5.02945848598297e-13\\
}; \addlegendentry{$\mathrm{M_6L_{20}}$}

\addplot [color=colorclassyblue, line width=0.5pt,loosely
dashdotted,mark=o,mark size=2,mark options={solid},mark repeat=1]
  table[row sep=crcr]{%
0.3125	1.06190229213169\\
0.217391304347826	1.17891888465003\\
0.15625	1.31853360417826\\
0.108695652173913	1.52320251357044\\
0.078125	1.52102605228755\\
0.0549450549450549	1.64780623561003\\
0.0390625	1.68710990279176\\
0.0274725274725275	1.45848815006373\\
0.01953125	1.66949530206489\\
0.0137741046831956	1.60007471533279\\
0.009765625	1.60523926363472\\
0.00689655172413793	1.2824715105952\\
0.0048828125	1.10490165519397\\
0.00345065562456867	1.38999247862795\\
0.00244140625	0.116021821690989\\
0.00172592336900242	0.000483456708317707\\
0.001220703125	1.45660552608923e-07\\
}; \addlegendentry{$\mathrm{M_6L_{10}}$}

\addplot [color=colorclassyorange, line width=0.5pt,mark=asterisk,mark
size=2,mark options={solid},mark repeat=1]
  table[row sep=crcr]{%
0.3125	0.269845085562075\\
0.217391304347826	0.0287740902541452\\
0.15625	0.00598535902887919\\
0.108695652173913	0.000464899349721709\\
0.078125	6.49508669185554e-05\\
0.0549450549450549	7.97747492861161e-06\\
0.0390625	1.03631923419894e-06\\
0.0274725274725275	1.25805095259265e-07\\
0.01953125	1.62687618157949e-08\\
0.0137741046831956	2.00462025471802e-09\\
0.009765625	2.54788425922393e-10\\
0.00689655172413793	3.72330841306388e-11\\
0.0048828125	4.31720502106715e-12\\
0.00345065562456867	8.12426190960238e-12\\
0.00244140625	7.03046271063629e-13\\
0.00172592336900242	1.94190526176589e-12\\
0.001220703125	6.61014646447738e-13\\
}; \addlegendentry{$\mathrm{S_6}$}

\addplot [color=colorclassyorange, line width=0.5pt,densely
dotted,mark=asterisk,mark size=2,mark options={solid},mark repeat=1]
  table[row sep=crcr]{%
0.3125	1.47678877565314\\
0.217391304347826	1.17008339090405\\
0.15625	1.56129895934373\\
0.108695652173913	1.57055480110608\\
0.078125	1.02672618053396\\
0.0549450549450549	0.706669933422044\\
0.0390625	0.43104684572394\\
0.0274725274725275	0.0355942849654915\\
0.01953125	9.17312698365334e-05\\
0.0137741046831956	2.7601968133944e-09\\
0.009765625	2.54670576935597e-10\\
0.00689655172413793	3.70961564548856e-11\\
0.0048828125	4.17594986944619e-12\\
0.00345065562456867	8.28837264972142e-12\\
0.00244140625	5.61370418612823e-13\\
0.00172592336900242	1.77607788228733e-12\\
0.001220703125	5.26840019969519e-13\\
}; \addlegendentry{$\mathrm{S_6L_{50}}$}

\addplot [color=colorclassyorange, line width=0.5pt,densely
dashed,mark=asterisk,mark size=2,mark options={solid},mark repeat=1]
  table[row sep=crcr]{%
0.3125	1.28459047802114\\
0.217391304347826	1.25607948289865\\
0.15625	1.62467131623098\\
0.108695652173913	1.54749043879984\\
0.078125	1.63457113943695\\
0.0549450549450549	1.56280380791037\\
0.0390625	1.54442973700244\\
0.0274725274725275	1.54023931120029\\
0.01953125	1.25212759737986\\
0.0137741046831956	1.78585717057363\\
0.009765625	0.770657430406791\\
0.00689655172413793	0.215609494943496\\
0.0048828125	8.74936250585927e-05\\
0.00345065562456867	1.0136776728465e-09\\
0.00244140625	5.58936314304851e-13\\
0.00172592336900242	1.77822507009347e-12\\
0.001220703125	5.2701834790726e-13\\
}; \addlegendentry{$\mathrm{S_6L_{20}}$}

\addplot [color=colorclassyorange, line width=0.5pt,loosely
dashdotted,mark=asterisk,mark size=1,mark options={solid},mark repeat=1]
  table[row sep=crcr]{%
0.3125	1.29568631346301\\
0.217391304347826	1.34086128188486\\
0.15625	1.44371643204014\\
0.108695652173913	1.10034435829984\\
0.078125	1.63752265755996\\
0.0549450549450549	1.67913954108928\\
0.0390625	1.50969114035941\\
0.0274725274725275	1.6867799566548\\
0.01953125	1.6592629550719\\
0.0137741046831956	1.62420244879812\\
0.009765625	1.38246946913015\\
0.00689655172413793	1.28993934138932\\
0.0048828125	1.10333745339554\\
0.00345065562456867	1.40038658668986\\
0.00244140625	0.11578380455519\\
0.00172592336900242	0.000482276351657018\\
0.001220703125	1.44655225957412e-07\\
}; \addlegendentry{$\mathrm{S_6L_{10}}$}

\end{axis}
\end{tikzpicture}%
    \caption{{\bf[Low oscillation regime ($V_\mathrm{S}$)]}: When applied to the low oscillatory regime of $V_{\mathrm{S}}$, the order four and order six Magnus expansions with simplified commutators,
    $\mathrm{S_4}$ and $\mathrm{S_6}$, have a similar error as the standard Magnus expansions, $\mathrm{M_4}$ and $\mathrm{M_6}$.
    Not much difference is made in this case ($V_{\mathrm{S}}$) by considering higher-order quadrature. On the left we use a Strang splitting (exponential midpoint rule)
    with $1024$ grid
    points and $5\times 10^7$ time steps as a reference. As we can see, the errors saturate around $10^{-8}$, which is the accuracy of this reference solution.
    For the other two plots we use $\mathrm{M_6}$ with $180$ grid points and $2\times 10^5$ time steps.}
    \label{fig:S6M6VS}
\end{figure}

{\bf Magnus--Lanczos}. In Figure~\ref{fig:S6M6VS} (right) and
Figure~\ref{fig:S6M6VE} (right) we show the errors in exponentiating the
order six Magnus expansions (with and without simplification of commutators) via 50, 20 and 10
Lanczos iterations respectively. The Magnus--Lanczos schemes with $n$ Lanczos
iterations are denoted with the postfix $\mathrm{L}_n$. It is evident from
these figures that exponentiation of Magnus expansions via Lanczos iterations
requires either a larger number of Lanczos iterations or smaller time steps
before we achieve the accuracy inherent in the Magnus expansion (i.e. the
accuracy of brute force exponentiation, $\mathrm{M_6}$ and $\mathrm{S_6}$).

Figures~\ref{fig:S6M6VS} (right) and \ref{fig:S6M6VE} (right) suggest that
there is scope for improvement in the efficient exponentiation of Magnus
expansions, particularly when it comes to large time steps, which can be
crucial for long term integration. In particular, it is worth exploring
whether Zassenhaus splittings confer an advantage here.

The convergence of Lanczos approximation to the exponential can be very
sensitive to the degree of spatial discretisation. $\mathrm{M_6L_{50}H}$, in
Figure~\ref{fig:S6M6VS} was run using 1024 grid points. Not only are the
Lanczos iterations more expensive in this case, but the convergence also
occurs much later. In general we need more iterations since the spectral
radius of the Magnus expansion is larger (growing quadratically with finer
spatial resolution).
%

\begin{figure}[tbh]
%
%
\definecolor{mycolor1}{rgb}{0.00000,0.44700,0.74100}%
\definecolor{mycolor2}{rgb}{0.85000,0.32500,0.09800}%
\definecolor{mycolor3}{rgb}{0.92900,0.69400,0.12500}%
\definecolor{mycolor4}{rgb}{0.49400,0.18400,0.55600}%
\definecolor{mycolor5}{rgb}{0.46600,0.67400,0.18800}%
\definecolor{mycolor6}{rgb}{0.63500,0.07800,0.18400}%
\begin{tikzpicture}

\begin{axis}[%
width=1.3in, height=2in, scale only axis, xmode=log, xmin=0.001, xmax=0.2,
xminorticks=true, ymode=log, ymin=1e-13, ymax=10, yminorticks=true, axis
background/.style={fill=white}, xticklabel style = {xshift=1.5ex},
xlabel={time step}, ylabel={$\mathrm{L}^2$ error}, ylabel shift=-8pt,
title={\texttt{expm} vs Strang 1e8},legend style={legend cell align=left,
align=left, draw=black,font=\scriptsize}, legend pos=south east, grid=both]

\addplot [color=black, densely dotted]
  table[row sep=crcr]{%
0.15625	596.046447753906\\
0.108695652173913	139.588194724862\\
0.078125	37.2529029846191\\
0.0549450549450549	9.11411382355726\\
0.0390625	2.3283064365387\\
0.0274725274725275	0.569632113972329\\
0.01953125	0.145519152283669\\
0.0137741046831956	0.035995939830998\\
0.009765625	0.00909494701772928\\
0.00689655172413793	0.00226218433691842\\
0.0048828125	0.00056843418860808\\
0.00345065562456867	0.000141777226236882\\
0.00244140625	3.5527136788005e-05\\
0.00172592336900242	8.87331780757525e-06\\
0.001220703125	2.22044604925031e-06\\
};\addlegendentry{$O(h^4)$}

\addplot [color=black, dashed]
  table[row sep=crcr]{%
0.15625	1455.19152283669\\
0.108695652173913	164.919889797804\\
0.078125	22.7373675443232\\
0.0549450549450549	2.7515136528068\\
0.0390625	0.35527136788005\\
0.0274725274725275	0.0429924008251063\\
0.01953125	0.00555111512312578\\
0.0137741046831956	0.000682936423418975\\
0.009765625	8.67361737988404e-05\\
0.00689655172413793	1.07594974407535e-05\\
0.0048828125	1.35525271560688e-06\\
0.00345065562456867	1.68814486939283e-07\\
0.00244140625	2.11758236813575e-08\\
0.00172592336900242	2.64319409124601e-09\\
0.001220703125	3.30872245021211e-10\\
}; \addlegendentry{$O(h^6)$}

\addplot [color=colorclassyblue,line width=0.5pt, densely dotted,
mark=square,mark size=2.5,mark options={solid},mark repeat=1]
  table[row sep=crcr]{%
0.15625	1.46444791160985\\
0.108695652173913	1.43871083967161\\
0.078125	1.58924009496002\\
0.0549450549450549	1.58425287912924\\
0.0390625	1.66589361870871\\
0.0274725274725275	0.560940730425351\\
0.01953125	0.0936287941291683\\
0.0137741046831956	0.0899602860671312\\
0.009765625	0.023584527096389\\
0.00689655172413793	0.0061488152814297\\
0.0048828125	0.00158081706621822\\
0.00345065562456867	0.000399007245013631\\
0.00244140625	0.000100582649112478\\
0.00172592336900242	2.51970856152281e-05\\
0.001220703125	6.31480147657455e-06\\
}; \addlegendentry{M4}

\addplot [color=colorclassyorange, line width=0.5pt, densely dotted,
mark=triangle*,mark size=2,mark options={solid},mark repeat=1]
  table[row sep=crcr]{%
0.15625	0.93625590797556\\
0.108695652173913	1.62267732757596\\
0.078125	1.54049008963352\\
0.0549450549450549	1.98339653697198\\
0.0390625	1.68490144794472\\
0.0274725274725275	0.830179324355683\\
0.01953125	0.311358922663057\\
0.0137741046831956	0.0957407549671055\\
0.009765625	0.0272675344472028\\
0.00689655172413793	0.00722920551483456\\
0.0048828125	0.00187672791622258\\
0.00345065562456867	0.00047593857418034\\
0.00244140625	0.000120262519417494\\
0.00172592336900242	3.01634908105291e-05\\
0.001220703125	7.56402249748967e-06\\
}; \addlegendentry{$\mathrm{S_4G_{11}}$}

\addplot [color=colorclassyblue,line width=0.5pt,mark=o,mark size=2,mark
options={solid},mark repeat=1]
  table[row sep=crcr]{%
0.15625	1.34109758562385\\
0.108695652173913	1.41488549847693\\
0.078125	0.563846894658796\\
0.0549450549450549	1.9560706525367\\
0.0390625	1.4537867601641\\
0.0274725274725275	0.115101311103859\\
0.01953125	0.00825274043427884\\
0.0137741046831956	0.00402191244478394\\
0.009765625	0.000498281307433897\\
0.00689655172413793	6.5667203732798e-05\\
0.0048828125	8.45709927938149e-06\\
0.00345065562456867	1.06876720554768e-06\\
0.00244140625	1.35585486224097e-07\\
0.00172592336900242	2.12821165944763e-08\\
0.001220703125	1.31843567399927e-08\\
}; \addlegendentry{$\mathrm{M_6}$}

\addplot [color=colorclassyorange, line
width=0.5pt,mark=asterisk,mark size=2,mark options={solid},mark repeat=1]
  table[row sep=crcr]{%
0.15625	1.33786227227718\\
0.108695652173913	1.45083319862689\\
0.078125	0.45550037851759\\
0.0549450549450549	1.38745830560852\\
0.0390625	0.993573103633415\\
0.0274725274725275	0.0400006311666004\\
0.01953125	0.0318811955104078\\
0.0137741046831956	0.000635883002798513\\
0.009765625	7.95216136003346e-05\\
0.00689655172413793	1.04815532015101e-05\\
0.0048828125	1.43635314883355e-06\\
0.00345065562456867	1.83880639683398e-07\\
0.00244140625	2.67404950330634e-08\\
0.00172592336900242	1.33687038261246e-08\\
0.001220703125	1.30345198300864e-08\\
}; \addlegendentry{$\mathrm{S_6}$}

\addplot [color=colorclassyorange, line
width=0.5pt, mark=triangle*,mark size=2,mark options={solid},mark repeat=1]
  table[row sep=crcr]{%
0.15625	1.90538136895669\\
0.108695652173913	0.454573666433515\\
0.078125	0.0485353601869392\\
0.0549450549450549	0.0355140992852013\\
0.0390625	0.00883765187502355\\
0.0274725274725275	0.00175609907117503\\
0.01953125	0.000307699137926943\\
0.0137741046831956	4.5085849755862e-05\\
0.009765625	6.3068979744476e-06\\
0.00689655172413793	8.23194489227567e-07\\
0.0048828125	1.07248203631513e-07\\
0.00345065562456867	1.87690831707442e-08\\
0.00244140625	1.31482424099866e-08\\
0.00172592336900242	1.30302006749071e-08\\
0.001220703125	1.30272413425749e-08\\
}; \addlegendentry{$\mathrm{S_6G_{11}}$}

\end{axis}
\end{tikzpicture}%
        \hspace{-0.5cm}
%
%
\definecolor{mycolor1}{rgb}{0.00000,0.44700,0.74100}%
\definecolor{mycolor2}{rgb}{0.85000,0.32500,0.09800}%
\definecolor{mycolor3}{rgb}{0.92900,0.69400,0.12500}%
\definecolor{mycolor4}{rgb}{0.49400,0.18400,0.55600}%
\definecolor{mycolor5}{rgb}{0.46600,0.67400,0.18800}%
\definecolor{mycolor6}{rgb}{0.63500,0.07800,0.18400}%
\begin{tikzpicture}

\begin{axis}[%
width=1.3in, height=2in, scale only axis, xmode=log, xmin=0.001, xmax=0.2,
xminorticks=true, ymode=log, ymin=1e-13, ymax=10, yminorticks=true, axis
background/.style={fill=white}, yticklabels={},xticklabel style =
{xshift=1.5ex}, xlabel={time step}, title={\texttt{expm} vs $\mathrm{M_6}$
2e5},legend style={legend cell align=left, align=left,
draw=black,font=\scriptsize}, legend pos=south east, grid=both]

\addplot [color=black, densely dotted]
  table[row sep=crcr]{%
0.15625	596.046447753906\\
0.108695652173913	139.588194724862\\
0.078125	37.2529029846191\\
0.0549450549450549	9.11411382355726\\
0.0390625	2.3283064365387\\
0.0274725274725275	0.569632113972329\\
0.01953125	0.145519152283669\\
0.0137741046831956	0.035995939830998\\
0.009765625	0.00909494701772928\\
0.00689655172413793	0.00226218433691842\\
0.0048828125	0.00056843418860808\\
0.00345065562456867	0.000141777226236882\\
0.00244140625	3.5527136788005e-05\\
0.00172592336900242	8.87331780757525e-06\\
0.001220703125	2.22044604925031e-06\\
};\addlegendentry{$O(h^4)$}

\addplot [color=black, dashed]
  table[row sep=crcr]{%
0.15625	1455.19152283669\\
0.108695652173913	164.919889797804\\
0.078125	22.7373675443232\\
0.0549450549450549	2.7515136528068\\
0.0390625	0.35527136788005\\
0.0274725274725275	0.0429924008251063\\
0.01953125	0.00555111512312578\\
0.0137741046831956	0.000682936423418975\\
0.009765625	8.67361737988404e-05\\
0.00689655172413793	1.07594974407535e-05\\
0.0048828125	1.35525271560688e-06\\
0.00345065562456867	1.68814486939283e-07\\
0.00244140625	2.11758236813575e-08\\
0.00172592336900242	2.64319409124601e-09\\
0.001220703125	3.30872245021211e-10\\
}; \addlegendentry{$O(h^6)$}

\addplot [color=colorclassyblue,line width=0.5pt, densely dotted,
mark=square,mark size=2.5,mark options={solid},mark repeat=1]
  table[row sep=crcr]{%
0.15625	1.46444790212231\\
0.108695652173913	1.43871083028851\\
0.078125	1.58924008469439\\
0.0549450549450549	1.58425286877393\\
0.0390625	1.66589360788791\\
0.0274725274725275	0.560940726814588\\
0.01953125	0.0936287934737801\\
0.0137741046831956	0.0899602853815855\\
0.009765625	0.0235845268541689\\
0.00689655172413793	0.00614881515179848\\
0.0048828125	0.00158081696627126\\
0.00345065562456867	0.000399007152645586\\
0.00244140625	0.000100582558080745\\
0.00172592336900242	2.5196992563304e-05\\
0.001220703125	6.31469848512718e-06\\
}; \addlegendentry{M4}

\addplot [color=colorclassyorange, line width=0.5pt, densely dotted,
mark=triangle*,mark size=2,mark options={solid},mark repeat=1]
  table[row sep=crcr]{%
0.15625	0.936255902004845\\
0.108695652173913	1.62267731696274\\
0.078125	1.54049007956006\\
0.0549450549450549	1.98339652408474\\
0.0390625	1.68490143693519\\
0.0274725274725275	0.830179318874504\\
0.01953125	0.311358920549969\\
0.0137741046831956	0.0957407542559656\\
0.009765625	0.0272675341814656\\
0.00689655172413793	0.00722920537944826\\
0.0048828125	0.00187672781562495\\
0.00345065562456867	0.000475938482562946\\
0.00244140625	0.000120262429586818\\
0.00172592336900242	3.01633994586064e-05\\
0.001220703125	7.56392288818574e-06\\
}; \addlegendentry{$\mathrm{S_4G_{11}}$}

\addplot [color=colorclassyblue,line width=0.5pt,mark=o,mark size=2,mark
options={solid},mark repeat=1]
  table[row sep=crcr]{%
0.15625	1.34109757692798\\
0.108695652173913	1.41488548930497\\
0.078125	0.563846891182799\\
0.0549450549450549	1.95607063975297\\
0.0390625	1.45378675074438\\
0.0274725274725275	0.115101310429537\\
0.01953125	0.00825274048361218\\
0.0137741046831956	0.00402191252506388\\
0.009765625	0.000498281388828176\\
0.00689655172413793	6.56672905857241e-05\\
0.0048828125	8.45717851375683e-06\\
0.00345065562456867	1.06877714397687e-06\\
0.00244140625	1.35047586049597e-07\\
0.00172592336900242	1.69188645818179e-08\\
0.001220703125	2.12187095513405e-09\\
}; \addlegendentry{$\mathrm{M_6}$}

\addplot [color=colorclassyorange, line
width=0.5pt,mark=asterisk,mark size=2,mark options={solid},mark repeat=1]
  table[row sep=crcr]{%
0.15625	1.33786226359781\\
0.108695652173913	1.45083318920927\\
0.078125	0.455500375875118\\
0.0549450549450549	1.38745829666663\\
0.0390625	0.993573097228243\\
0.0274725274725275	0.0400006307814206\\
0.01953125	0.0318811952572281\\
0.0137741046831956	0.000635883112383678\\
0.009765625	7.95215026691827e-05\\
0.00689655172413793	1.04814580816681e-05\\
0.0048828125	1.43621143303894e-06\\
0.00345065562456867	1.83336023127614e-07\\
0.00244140625	2.32702846112406e-08\\
0.00172592336900242	2.92159317753255e-09\\
0.001220703125	3.66730123791027e-10\\
}; \addlegendentry{$\mathrm{S_6}$}

\addplot [color=colorclassyorange, line
width=0.5pt, mark=triangle*,mark size=2,mark options={solid},mark repeat=1]
  table[row sep=crcr]{%
0.15625	1.90538135658045\\
0.108695652173913	0.454573663700766\\
0.078125	0.0485353597574888\\
0.0549450549450549	0.0355140989510568\\
0.0390625	0.00883765172794176\\
0.0274725274725275	0.00175609897021309\\
0.01953125	0.000307699048210622\\
0.0137741046831956	4.50857592492753e-05\\
0.009765625	6.30679604132303e-06\\
0.00689655172413793	8.23002975852821e-07\\
0.0048828125	1.06365708296721e-07\\
0.00345065562456867	1.34239023248805e-08\\
0.00244140625	1.69492737859523e-09\\
0.00172592336900242	2.12253324314226e-10\\
0.001220703125	2.6622875325243e-11\\
}; \addlegendentry{$\mathrm{S_6G_{11}}$}

\end{axis}
\end{tikzpicture}%
        \hspace{-0.5cm}
%
%
\definecolor{mycolor1}{rgb}{0.00000,0.44700,0.74100}%
\definecolor{mycolor2}{rgb}{0.85000,0.32500,0.09800}%
\definecolor{mycolor3}{rgb}{0.92900,0.69400,0.12500}%
\definecolor{mycolor4}{rgb}{0.49400,0.18400,0.55600}%
\definecolor{mycolor5}{rgb}{0.46600,0.67400,0.18800}%
\definecolor{mycolor6}{rgb}{0.30100,0.74500,0.93300}%
\begin{tikzpicture}

\begin{axis}[%
width=1.9in, height=2in, scale only axis, xmode=log, xmin=0.001, xmax=0.2,
xminorticks=true, ymode=log, ymin=1e-13, ymax=10, yminorticks=true, axis
background/.style={fill=white}, yticklabels={},xticklabel style =
{xshift=1.5ex}, xlabel={time step}, title={Lanczos vs $\mathrm{M_6}$
2e5},legend style={legend cell align=left, align=left,
draw=black,font=\scriptsize}, legend pos=south east, grid=both]

\addplot [color=colorclassyblue,line width=0.5pt,mark=o,mark size=2,mark
options={solid},mark repeat=1]
  table[row sep=crcr]{%
0.15625	1.34109757692798\\
0.108695652173913	1.41488548930497\\
0.078125	0.563846891182799\\
0.0549450549450549	1.95607063975297\\
0.0390625	1.45378675074438\\
0.0274725274725275	0.115101310429537\\
0.01953125	0.00825274048361218\\
0.0137741046831956	0.00402191252506388\\
0.009765625	0.000498281388828176\\
0.00689655172413793	6.56672905857241e-05\\
0.0048828125	8.45717851375683e-06\\
0.00345065562456867	1.06877714397687e-06\\
0.00244140625	1.35047586049597e-07\\
0.00172592336900242	1.69188645818179e-08\\
0.001220703125	2.12187095513405e-09\\
}; \addlegendentry{$\mathrm{M_6}$}

\addplot [color=colorclassyblue, line width=0.5pt,densely dotted,mark=o,mark
size=2,mark options={solid},mark repeat=1]
  table[row sep=crcr]{%
0.15625	1.57545801473433\\
0.108695652173913	1.41666318007449\\
0.078125	1.44954372567517\\
0.0549450549450549	1.33453457167324\\
0.0390625	1.39104462798816\\
0.0274725274725275	1.42399933568029\\
0.01953125	1.41459136622967\\
0.0137741046831956	0.0058786020707105\\
0.009765625	0.000498281388599537\\
0.00689655172413793	6.56672905773003e-05\\
0.0048828125	8.45717846779051e-06\\
0.00345065562456867	1.06877713777991e-06\\
0.00244140625	1.35047570433786e-07\\
0.00172592336900242	1.69188229906192e-08\\
0.001220703125	2.12184945008057e-09\\
}; \addlegendentry{$\mathrm{M_6L_{50}}$}

\addplot [color=colorclassyblue, line width=0.5pt,densely dashed,mark=o,mark
size=2,mark options={solid},mark repeat=1]
  table[row sep=crcr]{%
0.15625	1.60134233835096\\
0.108695652173913	1.69197353864725\\
0.078125	1.40841833747876\\
0.0549450549450549	1.41424846266557\\
0.0390625	1.43832841203956\\
0.0274725274725275	1.64669100649598\\
0.01953125	1.27615248750971\\
0.0137741046831956	1.45638921988195\\
0.009765625	1.41254776029257\\
0.00689655172413793	1.41432726623016\\
0.0048828125	0.0413424891217345\\
0.00345065562456867	0.000412184204786658\\
0.00244140625	1.35047535866803e-07\\
0.00172592336900242	1.69188257154264e-08\\
0.001220703125	2.12184807816701e-09\\
}; \addlegendentry{$\mathrm{M_6L_{20}}$}

\addplot [color=colorclassyblue, line width=0.5pt,loosely dashdotted,mark=o,mark
size=2,mark options={solid},mark repeat=1]
  table[row sep=crcr]{%
0.15625	1.25657212760898\\
0.108695652173913	1.35539908630107\\
0.078125	1.44826555710291\\
0.0549450549450549	1.4333885737076\\
0.0390625	1.13122061418963\\
0.0274725274725275	1.38107527341409\\
0.01953125	1.58238290296059\\
0.0137741046831956	1.40216548662166\\
0.009765625	1.46476522479032\\
0.00689655172413793	1.33396928480414\\
0.0048828125	1.43947919525137\\
0.00345065562456867	1.41399079596284\\
0.00244140625	0.903197337478608\\
0.00172592336900242	0.0690095336638816\\
0.001220703125	0.000316174305714055\\
}; \addlegendentry{$\mathrm{M_6L_{10}}$}

\addplot [color=colorclassyorange, line width=0.5pt,mark=triangle*,mark
size=2,mark options={solid},mark repeat=1]
  table[row sep=crcr]{%
0.15625	1.90538135658045\\
0.108695652173913	0.454573663700766\\
0.078125	0.0485353597574888\\
0.0549450549450549	0.0355140989510568\\
0.0390625	0.00883765172794176\\
0.0274725274725275	0.00175609897021309\\
0.01953125	0.000307699048210622\\
0.0137741046831956	4.50857592492753e-05\\
0.009765625	6.30679604132303e-06\\
0.00689655172413793	8.23002975852821e-07\\
0.0048828125	1.06365708296721e-07\\
0.00345065562456867	1.34239023248805e-08\\
0.00244140625	1.69492737859523e-09\\
0.00172592336900242	2.12253324314226e-10\\
0.001220703125	2.6622875325243e-11\\
}; \addlegendentry{$\mathrm{S_6G_{11}}$}

\addplot [color=colorclassyorange, line width=0.5pt,densely
dotted,mark=triangle*,mark size=2,mark options={solid},mark repeat=1]
  table[row sep=crcr]{%
0.15625	1.41255335773635\\
0.108695652173913	1.41891108060622\\
0.078125	1.50368900705527\\
0.0549450549450549	1.4290440640426\\
0.0390625	1.39260131989143\\
0.0274725274725275	1.40199234731561\\
0.01953125	0.251674732856498\\
0.0137741046831956	0.000170208295020938\\
0.009765625	6.30679599737376e-06\\
0.00689655172413793	8.2300296361273e-07\\
0.0048828125	1.063656891407e-07\\
0.00345065562456867	1.34238932161897e-08\\
0.00244140625	1.69490759444637e-09\\
0.00172592336900242	2.12278993277293e-10\\
0.001220703125	2.66086380055013e-11\\
}; \addlegendentry{$\mathrm{S_6G_{11}L_{50}}$}

\addplot [color=colorclassyorange, line width=0.5pt,densely dashed,mark=triangle*,mark
size=2,mark options={solid},mark repeat=1]
  table[row sep=crcr]{%
0.15625	1.58947882784952\\
0.108695652173913	1.39989676879941\\
0.078125	1.41873017147\\
0.0549450549450549	1.41810772637202\\
0.0390625	1.31180547517932\\
0.0274725274725275	1.57153786496422\\
0.01953125	1.37779133400864\\
0.0137741046831956	1.33019430246045\\
0.009765625	1.41321321351068\\
0.00689655172413793	1.41432545037535\\
0.0048828125	0.0417516632281952\\
0.00345065562456867	0.000158191461994412\\
0.00244140625	1.6954641321667e-09\\
0.00172592336900242	2.12282360403297e-10\\
0.001220703125	2.66041030127785e-11\\
}; \addlegendentry{$\mathrm{S_6G_{11}L_{20}}$}

\addplot [color=colorclassyorange, line width=0.5pt,loosely dashdotted,mark=triangle*,mark
size=1,mark options={solid},mark repeat=1]
  table[row sep=crcr]{%
0.15625	1.69977484296516\\
0.108695652173913	0.729466648250884\\
0.078125	0.993872408551149\\
0.0549450549450549	1.81743740123386\\
0.0390625	1.42260632607183\\
0.0274725274725275	1.77869936536265\\
0.01953125	1.57995953798228\\
0.0137741046831956	1.36481503114633\\
0.009765625	1.56116988504966\\
0.00689655172413793	1.29766541901104\\
0.0048828125	1.38754223323398\\
0.00345065562456867	1.41182419306541\\
0.00244140625	0.912743550080561\\
0.00172592336900242	0.0612524322145888\\
0.001220703125	0.000314568606976482\\
}; \addlegendentry{$\mathrm{S_6G_{11}L_{10}}$}

\end{axis}
\end{tikzpicture}%
        \\
    \caption{{\bf[Highly oscillatory regime ($V_\mathrm{E}$)]}: In the highly oscillatory regime of $V_{\mathrm{E}}$, we also include the 11
    Gauss--Legendre quadrature knots versions of $\mathrm{S_4}$ and $\mathrm{S_6}$, which are labeled with the postfix $\mathrm{G_{11}}$.
    A significant difference is made in this case ($V_{\mathrm{E}}$) by considering higher-order quadrature. Analytic expressions for integrals could be
    beneficial in such cases.
    On the left we use a Strang splitting (exponential midpoint rule)
    with $512$ grid points and $10^8$ time steps as a reference.
    For the other two plots we use $\mathrm{M_6}$ with $180$ grid points and $2\times 10^5$ time steps for the reference solution.}
    \label{fig:S6M6VE}
\end{figure}
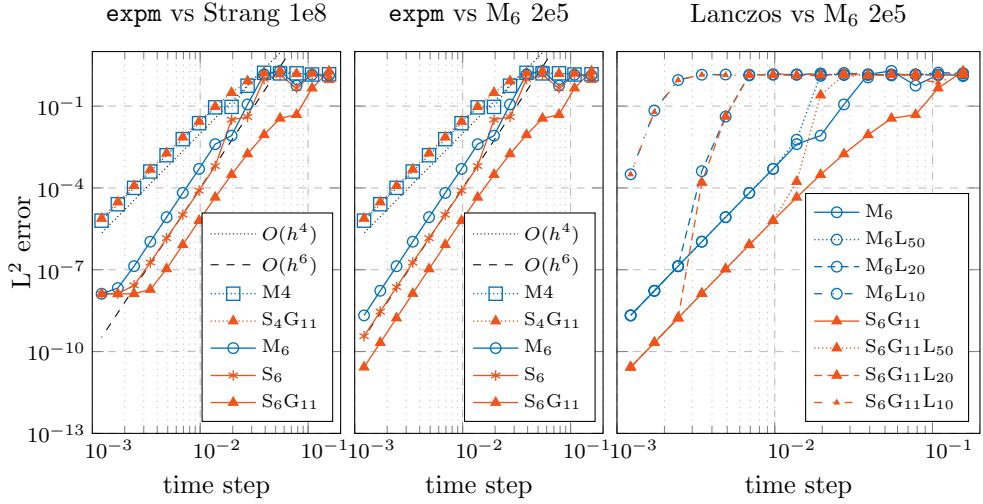

%


\begin{figure}[tbh]
%
%
\begin{tikzpicture}

\begin{axis}[%
width=1.85in, height=1.35in, scale only axis, xmode=log, xmin=0.001,
xmax=0.5, xminorticks=true, ymode=log, ymin=1, ymax=1e4, yminorticks=true,
axis background/.style={fill=white}, xlabel={time step},
ylabel={Computational Time (seconds)},ylabel style={font=\scriptsize},legend
style={legend cell align=left, align=left,title={$V_{\mathrm{S}}$},xticklabel
style = {xshift=1.5ex},ylabel style={font=\scriptsize},
draw=black,font=\scriptsize}, ytick={1e0,1e1,1e2,1e3,1e4,1e5}, yticklabel
style = {font=\scriptsize,xshift=0.25ex}, xticklabel style =
{font=\scriptsize,yshift=0.25ex}, legend pos=north east, grid=both] \addplot
[color=colorclassyblue,line width=0.5pt,mark=o,mark size=2,mark
options={solid},mark repeat=1]
  table[row sep=crcr]{%
0.3125	8\\
0.217391304347826	11\\
0.15625	15\\
0.108695652173913	20\\
0.078125	28\\
0.0549450549450549	39\\
0.0390625	54\\
0.0274725274725275	78\\
0.01953125	113\\
0.0137741046831956	156\\
0.009765625	218\\
0.00689655172413793	299\\
0.0048828125	437\\
0.00345065562456867	604\\
0.00244140625	852\\
0.00172592336900242	1203\\
0.001220703125	1696\\
}; \addlegendentry{$\mathrm{M_6L_{50}}$}

\addplot [color=colorclassyorange, line width=0.5pt,mark=triangle*,mark
size=2,mark options={solid},mark repeat=1]
  table[row sep=crcr]{%
0.3125	1\\
0.217391304347826	1\\
0.15625	1\\
0.108695652173913	2\\
0.078125	3\\
0.0549450549450549	4\\
0.0390625	5\\
0.0274725274725275	7\\
0.01953125	10\\
0.0137741046831956	15\\
0.009765625	20\\
0.00689655172413793	28\\
0.0048828125	39\\
0.00345065562456867	59\\
0.00244140625	87\\
0.00172592336900242	114\\
0.001220703125	156\\
}; \addlegendentry{$\mathrm{S_6L_{50}}$}

\end{axis}
\end{tikzpicture}%
        \hspace{0.5cm}
%
%
\begin{tikzpicture}

\begin{axis}[%
width=1.85in, height=1.35in, scale only axis, xmode=log, xmin=0.001,
xmax=0.2, xminorticks=true, ymode=log, ymin=1, ymax=1e4, yminorticks=true,
axis background/.style={fill=white}, xlabel={time step},
ylabel={Computational Time (seconds)}, ylabel
style={font=\scriptsize},xticklabel style =
{xshift=1.5ex},title={$V_{\mathrm{E}}$},legend style={legend cell align=left,
align=left, draw=black,font=\scriptsize}, yticklabel style =
{font=\scriptsize,xshift=0.25ex},xticklabel style =
{font=\scriptsize,yshift=0.25ex}, ytick={1e0,1e1,1e2,1e3,1e4,1e5}, legend
pos=north east, grid=both] \addplot [color=colorclassyblue,line
width=0.5pt,mark=o,mark size=2,mark options={solid},mark repeat=1]
  table[row sep=crcr]{%
0.15625	13\\
0.108695652173913	19\\
0.078125	26\\
0.0549450549450549	36\\
0.0390625	51\\
0.0274725274725275	72\\
0.01953125	101\\
0.0137741046831956	146\\
0.009765625	203\\
0.00689655172413793	294\\
0.0048828125	406\\
0.00345065562456867	577\\
0.00244140625	855\\
0.00172592336900242	1160\\
0.001220703125	1622\\
}; \addlegendentry{$\mathrm{M_6L_{50}}$}

\addplot [color=colorclassyorange, line width=0.5pt,mark=triangle*,mark
size=2,mark options={solid},mark repeat=1]
  table[row sep=crcr]{%
0.15625	2\\
0.108695652173913	3\\
0.078125	3\\
0.0549450549450549	5\\
0.0390625	7\\
0.0274725274725275	10\\
0.01953125	13\\
0.0137741046831956	19\\
0.009765625	27\\
0.00689655172413793	35\\
0.0048828125	47\\
0.00345065562456867	66\\
0.00244140625	93\\
0.00172592336900242	140\\
0.001220703125	191\\
}; \addlegendentry{$\mathrm{S_6G_{11}L_{50}}$}

\end{axis}
\end{tikzpicture}%
\\
    \caption{{\bf[Computational Time ($V_\mathrm{S}, V_\mathrm{E}$)]:} The absence of nested commutators
    in our proposed Magnus expansions results in a significant improvement in computational time
    of a corresponding Magnus--Lanczos scheme ($\mathrm{S_6L_{50}}$ is roughly 7 to 9 times faster here than $\mathrm{M_6L_{50}}$).
    This effect becomes more pronounced once we consider higher-order Magnus expansions.}
    \label{fig:CompTime}
\end{figure}
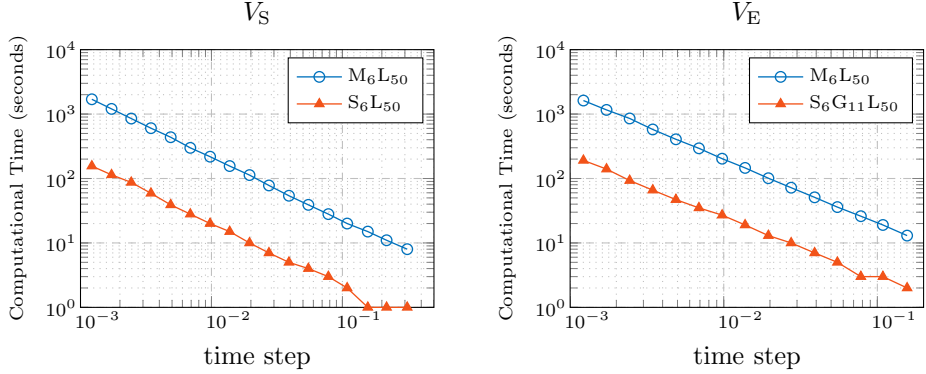

\begin{remark} Note that, for potentials of the form $V(x,t) =
V_0(x) + f(t)x$, the terms involving $\mom{1}{1},\mom{2}{1}$ and $\mom{3}{1}$
all vanish. This property, however, has not been exploited in the results
presented here.
\end{remark}

\section{Conclusions}

In this paper we have presented the derivation of integral-preserving Magnus--Lanczos methods
with simplified commutators (of arbitrarily high orders) for the
computation of the \schr equation featuring time-dependent potentials
\R{eq:schr} under the atomic scaling, $\hbar= 1$. In particular, we have
presented the 4th and 6th order methods (\R{eq:Omt3} and
(\ref{eq:Th4}--\ref{eq:Th4f}), respectively) and analysed their complexity in
terms of the number of FFTs required in each Lanczos iteration (see
\R{eq:Th4v}).

We find that the number of FFTs is much smaller than for standard Magnus
expansions where commutators appear explicitly -- our sixth-order method, for
instance, requires merely 8 FFTs for each Lanczos iteration. This speedup is
evident in numerical experiments (see Figure~\ref{fig:CompTime}, where we
find that our method is roughly 7 to 9 times faster than standard
Magnus--Lanczos methods). Moreover, the number of FFTs can be shown to grow
linearly with the order of the method we seek, so that the 8th order method
in this class of methods would require 10 FFTs and have a more pronounced
speedup over the standard Magnus expansion of order 8.

A concrete example of discretising the integrals via Gauss--Legendre
quadrature is also presented in (\ref{eq:FDu1}--\ref{eq:Lf}). However, as
stressed throughout, one of the major advantages of our approach is the
flexibility of choosing the method for approximating the integrals at the
very last stage. This is likely to prove highly beneficial in the case of
highly oscillatory potentials.

To illustrate this advantage, we present a numerical example featuring a
highly oscillatory potential, $V_{\mathrm{E}}$, in
Section~\ref{sec:experiments}. Here we find that our order-six Magnus
expansion using three Gauss--Legendres knots is roughly six times more
accurate than the standard Magnus expansion. This is improved significantly
by resorting to eleven Gauss--Legendre knots, resulting in the accuracy being
roughly 80 times higher than the standard Magnus expansion (see
Figure~\ref{fig:S6M6VE}).

\subsection{Future work}
{\bf Approximation of integrals.} It should be possible to improve upon the
accuracy and cost further by using analytic integrals, asymptotic
approximations or specialised highly-oscillatory quadrature.


{\bf Larger time steps.} We remind the reader that the accuracy inherent in
the Magnus expansion (as evident via direct exponentiation) is only reflected
in the Magnus--Lanczos methods when combined with a sufficient number of
Lanczos iterations (see Figure~\ref{fig:S6M6VS} (right) and
Figure~\ref{fig:S6M6VE} (right)). This is not much of a constraint when we
need high accuracy but can afford to work with moderately large to small time
steps since a reasonable number of Lanczos iterations suffices in this
regime. However, in applications where cost constraints trump accuracy
requirements and necessitate significantly larger time steps, the number of
Lanczos iterations required might become a concern.

In these regimes, it might be preferable to resort to Chebyshev expansions
which have been found to be more effective than Lanczos iterations for larger
time steps. Asymptotic splittings such as symmetric Zassenhaus splittings are
also likely to prove effective. This is because the number of Lanczos
iterations is dominated by the (spectral) size of the $\O{h}$ terms such as
$\ii \dt \dx^2$ and $-\ii \momO$ that arise from $\Theta^{[1]}(h)$, while the
cost of each Lanczos iteration in a high order expansion is dominated by the
trailing terms arising from $\Theta^{[k]}(h),k>1$. Symmetric Zassenhaus
splittings can effectively separate terms by powers of $h$ and are likely to
prove effective in decoupling these factors affecting the cost and accuracy
of Lanczos iterations.

{\bf Other equations.} It might be possible to extend some of these
techniques for other equations of quantum mechanics and certain linear
parabolic equations where Magnus expansions are employed. In particular, it
might be possible to simplify commutators\footnote{In applications where
there is no need to preserve skew-Hermiticity, it would suffice to expand in
terms of $f(x) \dx^k$ instead of working with $\ang{f}{k}$ (see Section~2).}
which could reduce cost, while preserving integrals in the case of highly
oscillatory forcing could increase accuracy.

\section*{Acknowledgments}
The work of KK in this project was financed by The National Center for
Science, based on decision no. 2016/22/M/ST1/00257.

\bibliographystyle{agsm}
\bibliography{Approach1}

\appendix
\section{Quadrature weights}
\label{sec:quadratureweights}
\newcommand{\knots}{K}
%
%
Once we have the values of $V$ at the set of knots $\knots$, the integrals
over the triangle can be approximated via \R{eq:Lf},
\[
\nonumber \Lf{f}{a}{b} \leadsto \sum_{j\in \knots} \sum_{k \in \knots} w_{jk}^f \left[
\K_{a} \MM{V}(\tau_j)\right] \left[\K_{b}
\MM{V}(\tau_k) \right] \ \text{or} \ \sum_{j\in \knots} \sum_{k \in \knots} w_{jk}^f \left[
 \MM{\dx^a V}(\tau_j)\right] \left[
\MM{\dx^b V}(\tau_k) \right],
\]
depending on whether exact derivatives $\dx^a V$ and $\dx^b V$ are available
or not (in the latter case we resort to numerical differentiation via
$\K_{a}$ and $\K_{b}$). As usual,  boldface denotes a vector of values
resulting from spatial discretisation.
 The weights required for three Gauss--Legendre quadrature knots
for the functions $\psi,\varphi_1,\varphi_2,\phi_1$ and $\phi_2$ are
\begin{eqnarray*}
w^{\psi} &=& \left(\frac{h}{6}\right)^3 \left\{\frac{1}{63} \left(
  \begin{array}{ccc}
    -139 & 26 & 239  \\
    26 & -304 & 26 \\
    239  & 26 & -139 \\
  \end{array}
\right)+5\sqrt{\frac{3}{5}}\left(
  \begin{array}{ccc}
    0 & 0 & -1 \\
    0 & 0 & 0 \\
    1 & 0 & 0 \\
  \end{array}
\right)\right\}\!,\\
w^{\varphi_1} &=& \frac{2}{7}\left(\frac{h}{6}\right)^4  \left\{ \left(
  \begin{array}{ccc}
    -11 & -62 & 136  \\
    190 & -128 & 190 \\
    136  & -62 & -11 \\
  \end{array}
\right) + \sqrt{\frac{3}{5}}
\left(
  \begin{array}{ccc}
    175 & 58 & -170  \\
    222 & 0 & -222 \\
    170  & -58 & -175 \\
  \end{array}
\right)\right\}\!,\\
w^{\varphi_2} &=& 2\left(\frac{h}{6}\right)^4  \left\{ \left(
  \begin{array}{ccc}
    -5 & -14 & 10  \\
    -2 & -32 & -2 \\
    10  & -14 & -5 \\
  \end{array}
\right) + \frac{1}{7}\sqrt{\frac{3}{5}}
\left(
  \begin{array}{ccc}
    145 & 134 & -90  \\
    -46 & 0 & 46 \\
    90  & -134 & -145 \\
  \end{array}
\right)\right\}\!,\\
w^{\phi_1} &=& \frac{2}{7}\left(\frac{h}{6}\right)^4  \left\{ \left(
  \begin{array}{ccc}
    -17 & 160 & -143  \\
    -92 & 184 & -92 \\
    -143  & 160 & -17 \\
  \end{array}
\right) + 6\sqrt{\frac{3}{5}}
\left(
  \begin{array}{ccc}
    25 & -34 & 30  \\
    -16 & 0 & 16 \\
    -30  & -34 & -25 \\
  \end{array}
\right)\right\}\!,\\
w^{\phi_2} &=& \left(\frac{h}{6}\right)^4  \left\{ 6 \left(
  \begin{array}{ccc}
    3 & 0 & -3  \\
    -4 & 8 & -4 \\
    -3  & 0 & 3 \\
  \end{array}
\right) + \frac{4}{7} \sqrt{\frac{3}{5}}
\left(
  \begin{array}{ccc}
    25 & -2 & 40  \\
    -48 & 0 & 48 \\
    -40  & 2 & -25 \\
  \end{array}
\right)\right\}\!.\\
\end{eqnarray*}


\end{document}